
\ifx\shlhetal\undefinedcontrolsequence\let\shlhetal\relax\fi

\input amstex
\expandafter\ifx\csname mathdefs.tex\endcsname\relax
  \expandafter\gdef\csname mathdefs.tex\endcsname{}
\else \message{Hey!  Apparently you were trying to
  \string\input{mathdefs.tex} twice.   This does not make sense.} 
\errmessage{Please edit your file (probably \jobname.tex) and remove
any duplicate ``\string\input'' lines}\endinput\fi




\catcode`\X=12\catcode`\@=11

\def\n@wcount{\alloc@0\count\countdef\insc@unt}
\def\n@wwrite{\alloc@7\write\chardef\sixt@@n}
\def\n@wread{\alloc@6\read\chardef\sixt@@n}
\def\r@s@t{\relax}\def\v@idline{\par}\def\@mputate#1/{#1}
\def\l@c@l#1X{\firstpart.#1}\def\gl@b@l#1X{#1}\def\t@d@l#1X{{}}

\def\crossrefs#1{\ifx\all#1\let\tr@ce=\all\else\def\tr@ce{#1,}\fi
   \n@wwrite\cit@tionsout\openout\cit@tionsout=\jobname.cit 
   \write\cit@tionsout{\tr@ce}\expandafter\setfl@gs\tr@ce,}
\def\setfl@gs#1,{\def\@{#1}\ifx\@\empty\let\next=\relax
   \else\let\next=\setfl@gs\expandafter\xdef
   \csname#1tr@cetrue\endcsname{}\fi\next}
\def\m@ketag#1#2{\expandafter\n@wcount\csname#2tagno\endcsname
     \csname#2tagno\endcsname=0\let\tail=\all\xdef\all{\tail#2,}
   \ifx#1\l@c@l\let\tail=\r@s@t\xdef\r@s@t{\csname#2tagno\endcsname=0\tail}\fi
   \expandafter\gdef\csname#2cite\endcsname##1{\expandafter
     \ifx\csname#2tag##1\endcsname\relax?\else\csname#2tag##1\endcsname\fi
     \expandafter\ifx\csname#2tr@cetrue\endcsname\relax\else
     \write\cit@tionsout{#2tag ##1 cited on page \folio.}\fi}
   \expandafter\gdef\csname#2page\endcsname##1{\expandafter
     \ifx\csname#2page##1\endcsname\relax?\else\csname#2page##1\endcsname\fi
     \expandafter\ifx\csname#2tr@cetrue\endcsname\relax\else
     \write\cit@tionsout{#2tag ##1 cited on page \folio.}\fi}
   \expandafter\gdef\csname#2tag\endcsname##1{\expandafter
      \ifx\csname#2check##1\endcsname\relax
      \expandafter\xdef\csname#2check##1\endcsname{}%
      \else\immediate\write16{Warning: #2tag ##1 used more than once.}\fi
      \multit@g{#1}{#2}##1/X%
      \write\t@gsout{#2tag ##1 assigned number \csname#2tag##1\endcsname\space
      on page \number\count0.}%
   \csname#2tag##1\endcsname}}

\def\multit@g#1#2#3/#4X{\def\t@mp{#4}\ifx\t@mp\empty%
      \global\advance\csname#2tagno\endcsname by 1 
      \expandafter\xdef\csname#2tag#3\endcsname
      {#1\number\csname#2tagno\endcsnameX}%
   \else\expandafter\ifx\csname#2last#3\endcsname\relax
      \expandafter\n@wcount\csname#2last#3\endcsname
      \global\advance\csname#2tagno\endcsname by 1 
      \expandafter\xdef\csname#2tag#3\endcsname
      {#1\number\csname#2tagno\endcsnameX}
      \write\t@gsout{#2tag #3 assigned number \csname#2tag#3\endcsname\space
      on page \number\count0.}\fi
   \global\advance\csname#2last#3\endcsname by 1
   \def\t@mp{\expandafter\xdef\csname#2tag#3/}%
   \expandafter\t@mp\@mputate#4\endcsname
   {\csname#2tag#3\endcsname\lastpart{\csname#2last#3\endcsname}}\fi}
\def\t@gs#1{\def\all{}\m@ketag#1e\m@ketag#1s\m@ketag\t@d@l p
\let\realscite\scite
\let\realstag\stag
   \m@ketag\gl@b@l r \n@wread\t@gsin
   \openin\t@gsin=\jobname.tgs \re@der \closein\t@gsin
   \n@wwrite\t@gsout\openout\t@gsout=\jobname.tgs }
\outer\def\localtags{\t@gs\l@c@l}
\outer\def\globaltags{\t@gs\gl@b@l}
\outer\def\newlocaltag#1{\m@ketag\l@c@l{#1}}
\outer\def\newglobaltag#1{\m@ketag\gl@b@l{#1}}

\newif\ifpr@ 
\def\m@kecs #1tag #2 assigned number #3 on page #4.%
   {\expandafter\gdef\csname#1tag#2\endcsname{#3}
   \expandafter\gdef\csname#1page#2\endcsname{#4}
   \ifpr@\expandafter\xdef\csname#1check#2\endcsname{}\fi}
\def\re@der{\ifeof\t@gsin\let\next=\relax\else
   \read\t@gsin to\t@gline\ifx\t@gline\v@idline\else
   \expandafter\m@kecs \t@gline\fi\let \next=\re@der\fi\next}
\def\pretags#1{\pr@true\pret@gs#1,,}
\def\pret@gs#1,{\def\@{#1}\ifx\@\empty\let\n@xtfile=\relax
   \else\let\n@xtfile=\pret@gs \openin\t@gsin=#1.tgs \message{#1} \re@der 
   \closein\t@gsin\fi \n@xtfile}

\newcount\sectno\sectno=0\newcount\subsectno\subsectno=0
\newif\ifultr@local \def\ultralocal{\ultr@localtrue}
\def\firstpart{\number\sectno}
\def\lastpart#1{\ifcase#1 \or a\or b\or c\or d\or e\or f\or g\or h\or 
   i\or k\or l\or m\or n\or o\or p\or q\or r\or s\or t\or u\or v\or w\or 
   x\or y\or z \fi}

\def\resetall{\global\advance\sectno by 1\subsectno=0
   \gdef\firstpart{\number\sectno}\r@s@t}
\def\resetsub{\global\advance\subsectno by 1
   \gdef\firstpart{\number\sectno.\number\subsectno}\r@s@t}
\def\newsection#1\par{\resetall\vskip0pt plus.3\vsize\penalty-250
   \vskip0pt plus-.3\vsize\bigskip\bigskip
   \message{#1}\leftline{\bf#1}\nobreak\bigskip}
\def\subsection#1\par{\ifultr@local\resetsub\fi
   \vskip0pt plus.2\vsize\penalty-250\vskip0pt plus-.2\vsize
   \bigskip\smallskip\message{#1}\leftline{\bf#1}\nobreak\medskip}


\newdimen\marginshift

\newdimen\margindelta
\newdimen\marginmax
\newdimen\marginmin

\def\margininit{       
\marginmax=3 true cm                  
				      
\margindelta=0.1 true cm              
\marginmin=0.1true cm                 
\marginshift=\marginmin
}    

\def\t@gsjj#1,{\def\@{#1}\ifx\@\empty\let\next=\relax\else\let\next=\t@gsjj
   \def\@@{p}\ifx\@\@@\else
   \expandafter\gdef\csname#1cite\endcsname##1{\citejj{##1}}
   \expandafter\gdef\csname#1page\endcsname##1{?}
   \expandafter\gdef\csname#1tag\endcsname##1{\tagjj{##1}}\fi\fi\next}
\newif\ifshowstuffinmargin
\showstuffinmarginfalse
\def\jjtags{\ifx\shlhetal\relax 
  \else
\ifx\shlhetal\undefinedcontrolseq
\else
\showstuffinmargintrue
\ifx\all\relax\else\expandafter\t@gsjj\all,\fi\fi \fi
}

\def\tagjj#1{\realstag{#1}\oldmginpar{\zeigen{#1}}}
\def\citejj#1{\rechnen{#1}\mginpar{\zeigen{#1}}}     

\def\rechnen#1{\expandafter\ifx\csname stag#1\endcsname\relax ??\else
                           \csname stag#1\endcsname\fi}

\newdimen\theight

\def\marginfont{\sevenrm}

\def\trymarginbox#1{\setbox0=\hbox{\marginfont\hskip\marginshift #1}%
		\global\marginshift\wd0 
		\global\advance\marginshift\margindelta}

\def \oldmginpar#1{%
\ifvmode\setbox0\hbox to \hsize{\hfill\rlap{\marginfont\quad#1}}%
\ht0 0cm
\dp0 0cm
\box0\vskip-\baselineskip
\else 
             \vadjust{\trymarginbox{#1}%
		\ifdim\marginshift>\marginmax \global\marginshift\marginmin
			\trymarginbox{#1}%
                \fi
             \theight=\ht0
             \advance\theight by \dp0    \advance\theight by \lineskip
             \kern -\theight \vbox to \theight{\rightline{\rlap{\box0}}%
\vss}}\fi}

\newdimen\upordown
\global\upordown=8pt
\font\tinyfont=cmtt8 
\def\mginpar#1{\smash{\hbox to 0cm{\kern-10pt\raise7pt\hbox{\tinyfont #1}\hss}}}
\def\mginpar#1{{\hbox to 0cm{\kern-10pt\raise\upordown\hbox{\tinyfont #1}\hss}}\global\upordown-\upordown}


\def\t@gsoff#1,{\def\@{#1}\ifx\@\empty\let\next=\relax\else\let\next=\t@gsoff
   \def\@@{p}\ifx\@\@@\else
   \expandafter\gdef\csname#1cite\endcsname##1{\zeigen{##1}}
   \expandafter\gdef\csname#1page\endcsname##1{?}
   \expandafter\gdef\csname#1tag\endcsname##1{\zeigen{##1}}\fi\fi\next}
\def\verbatimtags{\showstuffinmarginfalse
\ifx\all\relax\else\expandafter\t@gsoff\all,\fi}
\def\zeigen#1{\hbox{$\scriptstyle\langle$}#1\hbox{$\scriptstyle\rangle$}}


\def\margintag#1{\ifshowstuffinmargin\oldmginpar{\zeigen{#1}}\fi}

\def\(#1){\edef\dot@g{\ifmmode\ifinner(\hbox{\noexpand\etag{#1}})
   \else\noexpand\eqno(\hbox{\noexpand\etag{#1}})\fi
   \else(\noexpand\ecite{#1})\fi}\dot@g}

\newif\ifbr@ck
\def\eat#1{}
\def\[#1]{\br@cktrue[\br@cket#1'X]}
\def\br@cket#1'#2X{\def\temp{#2}\ifx\temp\empty\let\next\eat
   \else\let\next\br@cket\fi
   \ifbr@ck\br@ckfalse\br@ck@t#1,X\else\br@cktrue#1\fi\next#2X}
\def\br@ck@t#1,#2X{\def\temp{#2}\ifx\temp\empty\let\neext\eat
   \else\let\neext\br@ck@t\def\temp{,}\fi
   \def\teemp{#1}\ifx\teemp\empty\else\rcite{#1}\fi\temp\neext#2X}
\def\resetbr@cket{\gdef\[##1]{[\rtag{##1}]}}
\def\references{\resetbr@cket\newsection References\par}

\newtoks\symb@ls\newtoks\s@mb@ls\newtoks\p@gelist\n@wcount\ftn@mber
    \ftn@mber=1\newif\ifftn@mbers\ftn@mbersfalse\newif\ifbyp@ge\byp@gefalse
\def\defm@rk{\ifftn@mbers\n@mberm@rk\else\symb@lm@rk\fi}
\def\n@mberm@rk{\xdef\m@rk{{\the\ftn@mber}}%
    \global\advance\ftn@mber by 1 }
\def\rot@te#1{\let\temp=#1\global#1=\expandafter\r@t@te\the\temp,X}
\def\r@t@te#1,#2X{{#2#1}\xdef\m@rk{{#1}}}
\def\b@@st#1{{$^{#1}$}}\def\str@p#1{#1}
\def\symb@lm@rk{\ifbyp@ge\rot@te\p@gelist\ifnum\expandafter\str@p\m@rk=1 
    \s@mb@ls=\symb@ls\fi\write\f@nsout{\number\count0}\fi \rot@te\s@mb@ls}
\def\byp@ge{\byp@getrue\n@wwrite\f@nsin\openin\f@nsin=\jobname.fns 
    \n@wcount\currentp@ge\currentp@ge=0\p@gelist={0}
    \re@dfns\closein\f@nsin\rot@te\p@gelist
    \n@wread\f@nsout\openout\f@nsout=\jobname.fns }
\def\m@kelist#1X#2{{#1,#2}}
\def\re@dfns{\ifeof\f@nsin\let\next=\relax\else\read\f@nsin to \f@nline
    \ifx\f@nline\v@idline\else\let\t@mplist=\p@gelist
    \ifnum\currentp@ge=\f@nline
    \global\p@gelist=\expandafter\m@kelist\the\t@mplistX0
    \else\currentp@ge=\f@nline
    \global\p@gelist=\expandafter\m@kelist\the\t@mplistX1\fi\fi
    \let\next=\re@dfns\fi\next}
\def\symbols#1{\symb@ls={#1}\s@mb@ls=\symb@ls} 
\def\bigsymbol{\textstyle}
\symbols{\bigsymbol\ast,\dagger,\ddagger,\sharp,\flat,\natural,\star}
\def\ftnumbers{\ftn@mberstrue} \def\ftsymbols{\ftn@mbersfalse}
\def\paginal{\byp@ge} \def\resetftnumbers{\ftn@mber=1}
\def\ftnote#1{\defm@rk\expandafter\expandafter\expandafter\footnote
    \expandafter\b@@st\m@rk{#1}}

\long\def\jump#1\endjump{}
\def\ssum{\mathop{\lower .1em\hbox{$\textstyle\Sigma$}}\nolimits}

\def\qed{\nobreak\kern 1em \vrule height .5em width .5em depth 0em}
\def\newneq{\hbox{\rlap{\hbox to 1\wd9{\hss$=$\hss}}\raise .1em 
   \hbox to 1\wd9{\hss$\scriptscriptstyle/$\hss}}}
\def\subsetne{\setbox9 = \hbox{$\subset$}\mathrel{\hbox{\rlap
   {\lower .4em \newneq}\raise .13em \hbox{$\subset$}}}}
\def\supsetne{\setbox9 = \hbox{$\subset$}\mathrel{\hbox{\rlap
   {\lower .4em \newneq}\raise .13em \hbox{$\supset$}}}}

\def\vbar{\mathchoice{\vrule height6.3ptdepth-.5ptwidth.8pt\kern-.8pt}
   {\vrule height6.3ptdepth-.5ptwidth.8pt\kern-.8pt}
   {\vrule height4.1ptdepth-.35ptwidth.6pt\kern-.6pt}
   {\vrule height3.1ptdepth-.25ptwidth.5pt\kern-.5pt}}
\def\f@dge{\mathchoice{}{}{\mkern.5mu}{\mkern.8mu}}
\def\b@c#1#2{{\rm \mkern#2mu\vbar\mkern-#2mu#1}}
\def\b@b#1{{\rm I\mkern-3.5mu #1}}
\def\b@a#1#2{{\rm #1\mkern-#2mu\f@dge #1}}
\def\bb#1{{\count4=`#1 \advance\count4by-64 \ifcase\count4\or\b@a A{11.5}\or
   \b@b B\or\b@c C{5}\or\b@b D\or\b@b E\or\b@b F \or\b@c G{5}\or\b@b H\or
   \b@b I\or\b@c J{3}\or\b@b K\or\b@b L \or\b@b M\or\b@b N\or\b@c O{5} \or
   \b@b P\or\b@c Q{5}\or\b@b R\or\b@a S{8}\or\b@a T{10.5}\or\b@c U{5}\or
   \b@a V{12}\or\b@a W{16.5}\or\b@a X{11}\or\b@a Y{11.7}\or\b@a Z{7.5}\fi}}

\catcode`\X=11 \catcode`\@=12




\let\thischap\jobname

\def\partof#1{\csname returnthe#1part\endcsname}
\def\chapof#1{\csname returnthe#1chap\endcsname}

\def\setchapter#1,#2,#3;{%
  \expandafter\def\csname returnthe#1part\endcsname{#2}%
  \expandafter\def\csname returnthe#1chap\endcsname{#3}%
}

\setchapter 300a,A,II.A;
\setchapter 300b,A,II.B;
\setchapter 300c,A,II.C;
\setchapter 300d,A,II.D;
\setchapter 300e,A,II.E;
\setchapter 300f,A,II.F;
\setchapter 300g,A,II.G;
\setchapter  E53,B,N;
\setchapter  88r,B,I;
\setchapter  600,B,III;
\setchapter  705,B,IV;
\setchapter  734,B,V;

\def\cprefix#1{
\edef\theotherpart{\partof{#1}}\edef\theotherchap{\chapof{#1}}%
\ifx\theotherpart\thispart
   \ifx\theotherchap\thischap 
    \else 
     \theotherchap%
    \fi
   \else 
     \theotherchap\fi}

\def\sectioncite[#1]#2{%
     \cprefix{#2}#1}

\edef\thispart{\partof{\thischap}}
\edef\thischap{\chapof{\thischap}}

\def\lastpage of '#1' is #2.{\expandafter\def\csname lastpage#1\endcsname{#2}}


\def\spuriousreset{}


\expandafter\ifx\csname citeadd.tex\endcsname\relax
\expandafter\gdef\csname citeadd.tex\endcsname{}
\else \message{Hey!  Apparently you were trying to
\string\input{citeadd.tex} twice.   This does not make sense.} 
\errmessage{Please edit your file (probably \jobname.tex) and remove
any duplicate ``\string\input'' lines}\endinput\fi

\sectno=-1   
\localtags
\jjtags
\NoBlackBoxes
\define\mr{\medskip\roster}
\define\sn{\smallskip\noindent}
\define\mn{\medskip\noindent}
\define\bn{\bigskip\noindent}
\define\ub{\underbar}
\define\wilog{\text{without loss of generality}}
\define\ermn{\endroster\medskip\noindent}

\define \nl{\newline}
\magnification=\magstep 1
\documentstyle{amsppt}

{    
\catcode`@11

\ifx\alicetwothousandloaded@\relax
  \endinput\else\global\let\alicetwothousandloaded@\relax\fi

\gdef\subjclass{\let\savedef@\subjclass
 \def\subjclass##1\endsubjclass{\let\subjclass\savedef@
   \toks@{\def\usualspace{{\rm\enspace}}\eightpoint}%
   \toks@@{##1\unskip.}%
   \edef\thesubjclass@{\the\toks@
     \frills@{{\noexpand\rm2000 {\noexpand\it Mathematics Subject
       Classification}.\noexpand\enspace}}%
     \the\toks@@}}%
  \nofrillscheck\subjclass}
} 


\expandafter\ifx\csname alice2jlem.tex\endcsname\relax
  \expandafter\xdef\csname alice2jlem.tex\endcsname{\the\catcode`@}
\else \message{Hey!  Apparently you were trying to
\string\input{alice2jlem.tex}  twice.   This does not make sense.}
\errmessage{Please edit your file (probably \jobname.tex) and remove
any duplicate ``\string\input'' lines}\endinput\fi

\expandafter\ifx\csname bib4plain.tex\endcsname\relax
  \expandafter\gdef\csname bib4plain.tex\endcsname{}
\else \message{Hey!  Apparently you were trying to \string\input
  bib4plain.tex twice.   This does not make sense.}
\errmessage{Please edit your file (probably \jobname.tex) and remove
any duplicate ``\string\input'' lines}\endinput\fi

\def\renewcommand{\newcommand}	       
\edef\cite{\the\catcode`@}%
\catcode`@ = 11
\let\@oldatcatcode = \cite
\chardef\@letter = 11
\chardef\@other = 12
%
%
%
%
\def\@innerdef#1#2{\edef#1{\expandafter\noexpand\csname #2\endcsname}}%
%
%
\@innerdef\@innernewcount{newcount}%
\@innerdef\@innernewdimen{newdimen}%
\@innerdef\@innernewif{newif}%
\@innerdef\@innernewwrite{newwrite}%
%
%
%
\def\@gobble#1{}%
%
%
%
\ifx\inputlineno\@undefined
   \let\@linenumber = \empty 
\else
   \def\@linenumber{\the\inputlineno:\space}%
\fi
%
%
%
\def\@futurenonspacelet#1{\def\cs{#1}%
   \afterassignment\@stepone\let\@nexttoken=
}%
\begingroup 
\def\\{\global\let\@stoken= }%
\\ 
\endgroup
\def\@stepone{\expandafter\futurelet\cs\@steptwo}%
\def\@steptwo{\expandafter\ifx\cs\@stoken\let\@@next=\@stepthree
   \else\let\@@next=\@nexttoken\fi \@@next}%
\def\@stepthree{\afterassignment\@stepone\let\@@next= }%
%
%
%
\def\@getoptionalarg#1{%
   \let\@optionaltemp = #1%
   \let\@optionalnext = \relax
   \@futurenonspacelet\@optionalnext\@bracketcheck
}%
%
%
\def\@bracketcheck{%
   \ifx [\@optionalnext
      \expandafter\@@getoptionalarg
   \else
      \let\@optionalarg = \empty
      \expandafter\@optionaltemp
   \fi
}%
\def\@@getoptionalarg[#1]{%
   \def\@optionalarg{#1}%
   \@optionaltemp
}%
%
%
%
\def\@nnil{\@nil}%
\def\@fornoop#1\@@#2#3{}%
\def\@for#1:=#2\do#3{%
   \edef\@fortmp{#2}%
   \ifx\@fortmp\empty \else
      \expandafter\@forloop#2,\@nil,\@nil\@@#1{#3}%
   \fi
}%
\def\@forloop#1,#2,#3\@@#4#5{\def#4{#1}\ifx #4\@nnil \else
       #5\def#4{#2}\ifx #4\@nnil \else#5\@iforloop #3\@@#4{#5}\fi\fi
}%
\def\@iforloop#1,#2\@@#3#4{\def#3{#1}\ifx #3\@nnil
       \let\@nextwhile=\@fornoop \else
      #4\relax\let\@nextwhile=\@iforloop\fi\@nextwhile#2\@@#3{#4}%
}%
%
%
%
\@innernewif\if@fileexists
\def\@testfileexistence{\@getoptionalarg\@finishtestfileexistence}%
\def\@finishtestfileexistence#1{%
   \begingroup
      \def\extension{#1}%
      \immediate\openin0 =
         \ifx\@optionalarg\empty\jobname\else\@optionalarg\fi
         \ifx\extension\empty \else .#1\fi
         \space
      \ifeof 0
         \global\@fileexistsfalse
      \else
         \global\@fileexiststrue
      \fi
      \immediate\closein0
   \endgroup
}%
%
%
%
%
\def\bibliographystyle#1{%
   \@readauxfile
   \@writeaux{\string\bibstyle{#1}}%
}%
\let\bibstyle = \@gobble
%
%
\let\bblfilebasename = \jobname
\def\bibliography#1{%
   \@readauxfile
   \@writeaux{\string\bibdata{#1}}%
   \@testfileexistence[\bblfilebasename]{bbl}%
   \if@fileexists
      \nobreak
      \@readbblfile
   \fi
}%
\let\bibdata = \@gobble
%
%
\def\nocite#1{%
   \@readauxfile
   \@writeaux{\string\citation{#1}}%
}%
\@innernewif\if@notfirstcitation
%
%
\def\cite{\@getoptionalarg\@cite}%
%
%
\def\@cite#1{%
   \let\@citenotetext = \@optionalarg
   \printcitestart
   \nocite{#1}%
   \@notfirstcitationfalse
   \@for \@citation :=#1\do
   {%
      \expandafter\@onecitation\@citation\@@
   }%
   \ifx\empty\@citenotetext\else
      \printcitenote{\@citenotetext}%
   \fi
   \printcitefinish
}%
\newif\ifweareinprivate
\weareinprivatetrue
\ifx\shlhetal\undefinedcontrolseq\weareinprivatefalse\fi
\ifx\shlhetal\relax\weareinprivatefalse\fi
\def\@onecitation#1\@@{%
   \if@notfirstcitation
      \printbetweencitations
   \fi
   \expandafter \ifx \csname\@citelabel{#1}\endcsname \relax
      \if@citewarning
         \message{\@linenumber Undefined citation `#1'.}%
      \fi
     \ifweareinprivate
      \expandafter\gdef\csname\@citelabel{#1}\endcsname{%
\strut 
\vadjust{\vskip-\dp\strutbox
\vbox to 0pt{\vss\parindent0cm \leftskip=\hsize 
\advance\leftskip3mm
\advance\hsize 4cm\strut\openup-4pt 
\rightskip 0cm plus 1cm minus 0.5cm ?  #1 ?\strut}}
         {\tt
            \escapechar = -1
            \nobreak\hskip0pt\pfeilsw
            \expandafter\string\csname#1\endcsname
             \pfeilso
            \nobreak\hskip0pt
         }%
      }%
     \else  
      \expandafter\gdef\csname\@citelabel{#1}\endcsname{%
            {\tt\expandafter\string\csname#1\endcsname}
      }%
     \fi  
   \fi
   \csname\@citelabel{#1}\endcsname
   \@notfirstcitationtrue
}%
%
%
\def\@citelabel#1{b@#1}%
%
%
\def\@citedef#1#2{\expandafter\gdef\csname\@citelabel{#1}\endcsname{#2}}%
%
%
%
\def\@readbblfile{%
   \ifx\@itemnum\@undefined
      \@innernewcount\@itemnum
   \fi
   \begingroup
      \def\begin##1##2{%
         \setbox0 = \hbox{\biblabelcontents{##2}}%
         \biblabelwidth = \wd0
      }%
      \def\end##1{}
      %
      %
      \@itemnum = 0
      \def\bibitem{\@getoptionalarg\@bibitem}%
      \def\@bibitem{%
         \ifx\@optionalarg\empty
            \expandafter\@numberedbibitem
         \else
            \expandafter\@alphabibitem
         \fi
      }%
      \def\@alphabibitem##1{%
         \expandafter \xdef\csname\@citelabel{##1}\endcsname {\@optionalarg}%
         \ifx\biblabelprecontents\@undefined
            \let\biblabelprecontents = \relax
         \fi
         \ifx\biblabelpostcontents\@undefined
            \let\biblabelpostcontents = \hss
         \fi
         \@finishbibitem{##1}%
      }%
      \def\@numberedbibitem##1{%
         \advance\@itemnum by 1
         \expandafter \xdef\csname\@citelabel{##1}\endcsname{\number\@itemnum}%
         \ifx\biblabelprecontents\@undefined
            \let\biblabelprecontents = \hss
         \fi
         \ifx\biblabelpostcontents\@undefined
            \let\biblabelpostcontents = \relax
         \fi
         \@finishbibitem{##1}%
      }%
      \def\@finishbibitem##1{%
         \biblabelprint{\csname\@citelabel{##1}\endcsname}%
         \@writeaux{\string\@citedef{##1}{\csname\@citelabel{##1}\endcsname}}%
         \ignorespaces
      }%
      %
      %
      \let\em = \bblem
      \let\newblock = \bblnewblock
      \let\sc = \bblsc
      \frenchspacing
      \clubpenalty = 4000 \widowpenalty = 4000
      \tolerance = 10000 \hfuzz = .5pt
      \everypar = {\hangindent = \biblabelwidth
                      \advance\hangindent by \biblabelextraspace}%
      \bblrm
      \parskip = 1.5ex plus .5ex minus .5ex
      \biblabelextraspace = .5em
      \bblhook
      \input \bblfilebasename.bbl
   \endgroup
}%
%
%
\@innernewdimen\biblabelwidth
\@innernewdimen\biblabelextraspace
%
%
%
\def\biblabelprint#1{%
   \noindent
   \hbox to \biblabelwidth{%
      \biblabelprecontents
      \biblabelcontents{#1}%
      \biblabelpostcontents
   }%
   \kern\biblabelextraspace
}%
%
%
%
\def\biblabelcontents#1{{\bblrm [#1]}}%
%
%
\def\bblrm{\rm}%
%
%
\def\bblem{\it}%
%
%
\def\bblsc{\ifx\@scfont\@undefined
              \font\@scfont = cmcsc10
           \fi
           \@scfont
}%
%
%
\def\bblnewblock{\hskip .11em plus .33em minus .07em }%
%
%
\let\bblhook = \empty
%
%
%
\def\printcitestart{[}
\def\printcitefinish{]}
\def\printbetweencitations{, }
\def\printcitenote#1{, #1}
%
%
%
\let\citation = \@gobble
%
%
%
\@innernewcount\@numparams
%
%
\def\newcommand#1{%
   \def\@commandname{#1}%
   \@getoptionalarg\@continuenewcommand
}%
%
%
\def\@continuenewcommand{%
   \@numparams = \ifx\@optionalarg\empty 0\else\@optionalarg \fi \relax
   \@newcommand
}%
%
%
\def\@newcommand#1{%
   \def\@startdef{\expandafter\edef\@commandname}%
   \ifnum\@numparams=0
      \let\@paramdef = \empty
   \else
      \ifnum\@numparams>9
         \errmessage{\the\@numparams\space is too many parameters}%
      \else
         \ifnum\@numparams<0
            \errmessage{\the\@numparams\space is too few parameters}%
         \else
            \edef\@paramdef{%
               \ifcase\@numparams
                  \empty  No arguments.
               \or ####1%
               \or ####1####2%
               \or ####1####2####3%
               \or ####1####2####3####4%
               \or ####1####2####3####4####5%
               \or ####1####2####3####4####5####6%
               \or ####1####2####3####4####5####6####7%
               \or ####1####2####3####4####5####6####7####8%
               \or ####1####2####3####4####5####6####7####8####9%
               \fi
            }%
         \fi
      \fi
   \fi
   \expandafter\@startdef\@paramdef{#1}%
}%
%
%
%
%
\def\@readauxfile{%
   \if@auxfiledone \else 
      \global\@auxfiledonetrue
      \@testfileexistence{aux}%
      \if@fileexists
         \begingroup
            \endlinechar = -1
            \catcode`@ = 11
            \input \jobname.aux
         \endgroup
      \else
         \message{\@undefinedmessage}%
         \global\@citewarningfalse
      \fi
      \immediate\openout\@auxfile = \jobname.aux
   \fi
}%
%
%
\newif\if@auxfiledone
\ifx\noauxfile\@undefined \else \@auxfiledonetrue\fi
%
%
%
%
\@innernewwrite\@auxfile
\def\@writeaux#1{\ifx\noauxfile\@undefined \write\@auxfile{#1}\fi}%
%
%
%
\ifx\@undefinedmessage\@undefined
   \def\@undefinedmessage{No .aux file; I won't give you warnings about
                          undefined citations.}%
\fi
%
%
\@innernewif\if@citewarning
\ifx\noauxfile\@undefined \@citewarningtrue\fi
%
%
%
\catcode`@ = \@oldatcatcode

\def\pfeilso{\leavevmode
            \vrule width 1pt height9pt depth 0pt\relax
           \vrule width 1pt height8.7pt depth 0pt\relax
           \vrule width 1pt height8.3pt depth 0pt\relax
           \vrule width 1pt height8.0pt depth 0pt\relax
           \vrule width 1pt height7.7pt depth 0pt\relax
            \vrule width 1pt height7.3pt depth 0pt\relax
            \vrule width 1pt height7.0pt depth 0pt\relax
            \vrule width 1pt height6.7pt depth 0pt\relax
            \vrule width 1pt height6.3pt depth 0pt\relax
            \vrule width 1pt height6.0pt depth 0pt\relax
            \vrule width 1pt height5.7pt depth 0pt\relax
            \vrule width 1pt height5.3pt depth 0pt\relax
            \vrule width 1pt height5.0pt depth 0pt\relax
            \vrule width 1pt height4.7pt depth 0pt\relax
            \vrule width 1pt height4.3pt depth 0pt\relax
            \vrule width 1pt height4.0pt depth 0pt\relax
            \vrule width 1pt height3.7pt depth 0pt\relax
            \vrule width 1pt height3.3pt depth 0pt\relax
            \vrule width 1pt height3.0pt depth 0pt\relax
            \vrule width 1pt height2.7pt depth 0pt\relax
            \vrule width 1pt height2.3pt depth 0pt\relax
            \vrule width 1pt height2.0pt depth 0pt\relax
            \vrule width 1pt height1.7pt depth 0pt\relax
            \vrule width 1pt height1.3pt depth 0pt\relax
            \vrule width 1pt height1.0pt depth 0pt\relax
            \vrule width 1pt height0.7pt depth 0pt\relax
            \vrule width 1pt height0.3pt depth 0pt\relax}

\def\pfeilsw{ \leavevmode 
            \vrule width 1pt height0.3pt depth 0pt\relax
            \vrule width 1pt height0.7pt depth 0pt\relax
            \vrule width 1pt height1.0pt depth 0pt\relax
            \vrule width 1pt height1.3pt depth 0pt\relax
            \vrule width 1pt height1.7pt depth 0pt\relax
            \vrule width 1pt height2.0pt depth 0pt\relax
            \vrule width 1pt height2.3pt depth 0pt\relax
            \vrule width 1pt height2.7pt depth 0pt\relax
            \vrule width 1pt height3.0pt depth 0pt\relax
            \vrule width 1pt height3.3pt depth 0pt\relax
            \vrule width 1pt height3.7pt depth 0pt\relax
            \vrule width 1pt height4.0pt depth 0pt\relax
            \vrule width 1pt height4.3pt depth 0pt\relax
            \vrule width 1pt height4.7pt depth 0pt\relax
            \vrule width 1pt height5.0pt depth 0pt\relax
            \vrule width 1pt height5.3pt depth 0pt\relax
            \vrule width 1pt height5.7pt depth 0pt\relax
            \vrule width 1pt height6.0pt depth 0pt\relax
            \vrule width 1pt height6.3pt depth 0pt\relax
            \vrule width 1pt height6.7pt depth 0pt\relax
            \vrule width 1pt height7.0pt depth 0pt\relax
            \vrule width 1pt height7.3pt depth 0pt\relax
            \vrule width 1pt height7.7pt depth 0pt\relax
            \vrule width 1pt height8.0pt depth 0pt\relax
            \vrule width 1pt height8.3pt depth 0pt\relax
            \vrule width 1pt height8.7pt depth 0pt\relax
            \vrule width 1pt height9pt depth 0pt\relax
      }


\def\widestnumber#1#2{}

\def\citewarning#1{\ifx\shlhetal\relax 
    \else
    \par{#1}\par
    \fi
}

\def\rm{\fam0 \tenrm}

\def\fakesubhead#1\endsubhead{\bigskip\noindent{\bf#1}\par}



%
%
%

%

\font\textrsfs=rsfs10
\font\scriptrsfs=rsfs7
\font\scriptscriptrsfs=rsfs5

\newfam\rsfsfam
\textfont\rsfsfam=\textrsfs
\scriptfont\rsfsfam=\scriptrsfs
\scriptscriptfont\rsfsfam=\scriptscriptrsfs

\edef\oldcatcodeofat{\the\catcode`\@}
\catcode`\@11

\def\Cal@@#1{\noaccents@ \fam \rsfsfam #1}

\catcode`\@\oldcatcodeofat


\expandafter\ifx \csname margininit\endcsname \relax\else\margininit\fi

\long\def\red#1\endred{}
\long\def\green#1\endgreen{}
\long\def\blue#1\endblue{}
\long\def\private#1\endprivate{}

\def\endred{ \unmatched endred! }
\def\endgreen{ \unmatched endgreen! }
\def\endblue{ \unmatched endblue! }
\def\endprivate{ \unmatched endprivate! }

\ifx\latexcolors\undefinedcs\def\latexcolors{}\fi

\def\emptycs{}
\def\evaluatelatexcolors{%
        \ifx\latexcolors\emptycs\else
        \expandafter\xxevaluate\latexcolors\xxfertig\evaluatelatexcolors\fi}
\def\xxevaluate#1,#2\xxfertig{\setupthiscolor{#1}%
        \def\latexcolors{#2}}


\font\smallfont=cmsl7
\def\rutgerscolor{\ifmmode\else\endgraf\fi\smallfont
\advance\leftskip0.5cm\relax}
\def\setupthiscolor#1{\edef\tmptmpcs{\noexpand\bgroup\noexpand\rutgerscolor
\noexpand\def\noexpand\currentcolor{#1}%
\noexpand}%
\expandafter\let\csname#1\endcsname\tmptmpcs
\def\tmptmpcs{\checkColorUnmatched{#1}\popthecolor}
\expandafter\let\csname end#1\endcsname\tmptmpcs}

\def\checkColorUnmatched#1{\def\expectcolor{#1}%
    \ifx\expectcolor\currentcolor   
    \else \edef\failhere{\noexpand\tryingToClose '\currentcolor' with end\expectcolor}\failhere\fi}

\def\currentcolor{???}

\def\popthecolor{\ifmmode\else\endgraf\fi\egroup}

\expandafter\def\csname#1\endcsname{}

\evaluatelatexcolors

 \let\outerhead\head
 \def\head{\innerhead}
 \let\innerhead\outerhead

 \let\outersubhead\subhead
 \def\subhead{\innersubhead}
 \let\innersubhead\outersubhead

 \let\outersubsubhead\subsubhead
 \def\subsubhead{\innersubsubhead}
 \let\innersubsubhead\outersubsubhead

 \let\outerproclaim\proclaim
 \def\proclaim{\innerproclaim}
 \let\innerproclaim\outerproclaim

 %
 %
 %
 %

\def\demo#1{\medskip\noindent{\it #1.\/}}
\def\enddemo{\smallskip}

\def\remark#1{\medskip\noindent{\it #1.\/}}
\def\endremark{\smallskip}

\pageheight{8.5truein}
\topmatter
\title{What majority decisions are possible} \endtitle
\rightheadtext{Majority Decisions}
\author {Saharon Shelah \thanks {\null\newline 
Partially supported by the United States-Israel Binational Science
Foundation. Publication 816. \null\newline
I would like to thank Alice Leonhardt for the beautiful typing. 
} \endthanks} \endauthor 

 \affil{The Hebrew University of Jerusalem \\
Einstein Institute of Mathematics\\
Edmond J. Safra Campus, Givat Ram \\
 Jerusalem 91904, Israel
 \medskip
 Department of Mathematics \\
 Hill Center-Busch Campus \\
 Rutgers, The State University of New Jersey \\
 110 Frelinghuysen Road \\
 Piscataway, NJ 08854-8019  USA} \endaffil

\abstract  Suppose we are given a family of choice functions on pairs
from a given finite set (with at least three elements) closed under
permutations of the given set.  The set is considered the set of
alternatives (say candidates for an office).  The question is, what
are the choice functions $\bold c$ on pairs of this set of the
following form: for some (finite) family of ``voters", each having a
preference, i.e., a choice from each pair from the given family,
$\bold c\{x,y\}$ is chosen by the preference of the majority of
voters.  We give full characterization.  \endabstract
\endtopmatter
\document

\newpage

\head {Annotated Content} \endhead
 \spuriousreset
\bn
\S0 Introduction
\bn
\S1 Basic definitions and facts
\bn
\S2 When every majority choice is possible: a characterization
\bn
\S3 Balanced choice functions
\mr
\item "{${{}}$}" [We characterize what the majority choice can be on
pr-$c \ell({\Cal D})$ for ${\Cal D} \subseteq {\frak C}$ which is
balanced, i.e., does not fall under \S2.  We get the full answer.]
\endroster
\newpage

\head {\S0 Introduction} \endhead  \resetall \sectno=0
 \spuriousreset
\bigskip

Condorcet's ``paradox'' demonstrates that given three candidates A, B and C,
the majority rule may result in the 
society preferring A to B , B to C and C to A.
McGarvey \cite {McG53} proved a far-reaching extension of Condorcet's paradox:
For every asymmetric relation $R$ on a finite set $M$ of
candidates there is a strict-preferences
(linear orders, no ties) voter profile that has the relation $R$
as its strict simple majority relation. In other words,
for every assymetric relation (equivalently, a tournament)
$R$ on a set $M$ of $m$ elements
there are $n$ linear order relations on $M$, $R_1,R_2, \dots, R_n$
such that for every $a,b \in M$, $aRb$ if and only if
$$|\{ i: aR_ib\}|> n/2.$$ McGarvey's proof gave $n = m(m-1)$.
Stearns \cite {Ste59} found a construction with
$n = m$ and noticed that a simple 
counting argument implies that $n$ must be at least  $m/\log m$.
Erd\H{o}s and Moser \cite {ErMo64} were able to give a construction
with $n=O(m/\log m)$. Alon \cite {Alo02} showed that for some constant
$c_1>0$ we can find $R_1,\dots, R_n$
with $$|\{ i: aR_ib\}|> (1/2+ c_1/\sqrt {n})n,$$ and
that this is no longer the case if $c_1$ is
replaced with another constant $c_2 > c_1$.
\sn
Gil Kalai asked to what extent the assertion of McGarvey's theorem
holds if we replace the set of order relations
by an arbitrary isomorphism class of choice functions on pairs of
elements (see Definition \scite{0.4}).
Namely, the question is to characterize under which conditions clause
(A) of \scite{0.1} below holds (i.e., question \scite{1.3}).
\nl
Instead of choice functions we can speak on tournaments, see
observation \scite{0.5}. \nl

The main result is (follows from \scite{2.1})
\proclaim{\stag{0.1} Theorem}  Let $X$ be a finite set and ${\frak D}$
be a non-empty family of choice functions for $\binom X 2$ closed under
permutations of $X$.  \ub{Then} the following conditions are equivalent:
\mr
\item "{$(A)$}"  for any choice function $c$ on $\binom X2$ we can
find a finite set $J$ and $c_j \in {\frak D}$ for $j \in J$ such that
for any $x \ne y \in X$:
$$
c\{x,y\} = y \text{ \ub{iff} } |J|/2 < |\{j \in J:c_j\{x,y\} = y\}|
$$
(so equality never occurs)
\sn
\item "{$(B)$}"  for some $c \in {\frak D}$ and some $x \in X$ we have
$|\{y:c\{x,y\} = y\}| \ne (|X|-1)/2$.
\endroster
\endproclaim
\bn
Gil Kalai further asks \nl
\margintag{0.Q}\ub{\stag{0.Q} Question}:  1) In \scite{0.1} can we bound $|J|$? 
\nl
2) What is the result of demanding a ``non-trivial majority"? (say 51\%?)
\mn
Under \scite{0.1} it seems reasonable to characterize what can be
$\{c:c$ a choice function for pairs from $X$ gotten as in clause (A)
of \scite{0.1} using $c_j \in {\frak D}\}$, when we vary ${\frak D}$,
so \scite{0.1} tells us for which sets ${\frak D}$ the resulting family
is maximal.
\nl
We then give in \scite{b.5} a 
complete solution also to the question: what is the closure of
a set of choice functions by majority; in fact, there are just two.

We also may allow each ``voter" to abstain, this means
that his choice function is only partial.  We hope to deal with this
elsewhere, but there are more cases, e.g. of course, if 
all voters have no opinion on any pair, majority
discussion will always be a draw (giving a third possibility).  
Note that now we consider also
majority decisions which give a draw in some of the cases.  The
present work was
present in the conference in honour of Michael O. Rabin, Summer 2005.
\nl
I thank Gil for the stimulating discussion and writing the historical
background and the referee and Mor Doron for pointing out errors and
helping in proofreadings. 
\nl
An earlier version is \cite{Sh:E37}.
\bn
\margintag{0.2}\ub{\stag{0.2} Notation}:  Let $n,m,k,\ell,i,j$ denote natural numbers. \nl
Let $r,s,t,a,b$ denote real numbers.  \nl
Let $x,y,z,u,v,w$ denote members of the finite set $X$. \nl
Let $\binom Xk$ be the family of subsets of $X$ with exactly $k$
members. \nl
Let $c,d$ denote partial choice functions on $\binom X2$. \nl
Let conv$(A)$ be the convex hull of $A$, here for $A \subseteq \Bbb R
\times \Bbb R$. \nl
Let Per$(X)$ be the set of permutations of $X$.
\bn
The ``translation" to tournaments is not really used, still we explain it.
\definition{\stag{0.4} Definition}  1) We say that $c$ is a choice
[partial choice] function for $\binom Xk$ if $c$ is a function with
domain $\binom Xk$ [domain $\subseteq \binom Xk$] such that $x \in
\text{ Dom}(c) \Rightarrow c(x) \in x$.
\nl
2) We say $c$ is a choice function for pairs from $X$ \ub{when} $c$ is
a choice function for $\binom X2$.
\nl
3) If $c$ is a partial choice function for pairs from $X$, let
Tor$[c]$ be the following directed graph:
\mr
\item "{$(a)$}"  the set of nodes is $X$
\sn
\item "{$(b)$}"  the set of edges is $\{(x,y):x \ne y$ are from $X$
and $c\{x,y\}=y\}$.
\ermn
3A) Let $\bold G_c$ be the non-directed graph derived from Tor$[c]$.
\nl
4) Let $c_1,c_2$ be choice functions for pairs from $X$.  We say $\pi$
is an isomorphism from $c_1$ onto $c_2$ \ub{if} $\pi$ is a permutation
of $X$ such that for every $x,y \in X$ we have $c_1\{x,y\} = y
\Leftrightarrow c_2\{\pi(x),\pi(y)\} = \pi(y)$.
\enddefinition
\bigskip

\demo{\stag{0.5} Observation}  1) For any set $X$, the mapping $c
\mapsto \text{ Tor}[c]$ is a one-to-one mapping from the set of choice
functions for pairs from $X$ onto the set of tournaments on $X$.
\nl
1A)  It is also a one to one map from the set of partial choice
functions onto the set of directed graphs on $X$ (so for $x \ne y$
maybe $(x,y)$ is an edge maybe $(y,x)$ is an edge but not both and
maybe none).
\nl
2) For choice functions $c_1,c_2$ of pairs from $X$; we have $c_1,c_2$
are isomorphic iff Tor$[c_1]$, Tor$[c_2]$ are isomorphic tournaments.
\nl
2A) Similarly for partial choice functions and directed graphs.
\enddemo
\newpage

\head {\S1 Basic definitions and facts} \endhead  \resetall \sectno=1
 \spuriousreset
\bigskip

\demo{\stag{1.1} Hypothesis}  Assume
\mr
\item "{$(a)$}"  $X$ is a (fixed) finite set with $\bold n \ge 3$
members, i.e. $\bold n = |X|$.
\sn
\item "{$(b)$}"  ${\frak C} = {\frak C}^1 = {\frak C}^1_X$ is the set of
partial choice functions on $\binom X2$, see Definition \scite{0.4}(1);
when $c\{x,y\}$ is not defined it is interpreted as abstaining or
having no preference. \nl
Let ${\frak C}^0 = {\frak C}^{\text{full}} = {\frak C}^{\text{full}}_X
= {\frak C}^0_X$ be the set of $c \in {\frak C}^1_X$ which are full,
i.e., Dom$(c) = \binom X2 = \{\{x,y\}:x \ne y \in X\}$.
\sn
\item "{$(c)$}"  ${\Cal C},{\Cal D}$ vary on subsets of ${\frak C}$
\sn
\item "{$(d)$}"  ${\frak D}$ vary on non-empty subsets of 
${\frak C}$ which are symmetric where
\endroster
\enddemo
\bigskip

\definition{\stag{1.1A} Definition}  1) ${\Cal C} \subseteq {\frak
C}$ is symmetric \ub{if} it is closed under permutations of $X$ 
(i.e. for every $\pi \in \text{ Per}(X)$ the permutation $\hat \pi$ maps
${\Cal C}$ onto itself where $\pi$ induces $\hat \pi$, 
a permutation of ${\frak C}$, that is $c_1 = c^\pi_2$ or $c_1
= \hat \pi c_2$ mean that: $x_1 = \pi(x_2),y_1 = \pi(y_2)$ implies
$c_1\{x_1,y_1\} = y_1 \Leftrightarrow c_2\{x_2,y_2\} = y_2$).
\nl
2) For ${\Cal D} \subseteq {\frak C}$ and $x \ne y \in X$ let 
${\Cal D}_{x,y} = \{d \in {\Cal D}:d\{x,y\} = y \}$.
\enddefinition
\bigskip

\definition{\stag{1.2} Definition}  For ${\Cal D} \subseteq {\frak
C}$ let maj-$cl({\Cal D})$ be the set of $d \in {\frak C}$ such that
for some real numbers $r_c = r_c[d] \in [0,1]_{\Bbb R}$ for $c \in {\Cal D}$
satisfying $\underset{c \in {\Cal D}} {}\to \Sigma r_c = 1$ we have
\footnote{note that there is no a priori reason to assume that ${\Cal D}_2 =
\text{ maj}-c \ell({\Cal D}_1)$ implies ${\Cal D}_2 = \text{ maj}-c
\ell({\Cal D}_2)$} 

$$
\align
d\{x,y\} = &x \Leftrightarrow \frac 12 < \Sigma\{r_c:c\{x,y\} = x
\text{ and } c \in {\Cal D}\} \\
  &+\Sigma\{r_c/2:c\{x,y\} \text{ is undefined and } c \in {\Cal D}\}.
\endalign
$$
\enddefinition
\bigskip

\remark{Remark}  1) Clearly maj is for majority.  At first glance this
is not the same as the problem stated in the introduction but easily
they are equivalent (see clause (c) of \scite{2.1}).
\nl
2) Note that if we deal with full choice functions only, as
originally, then we require that the sum is never $\frac 12$.
\nl
3) Modulo the equivalence above, Kalai's original question was 
\endremark
\sn
\margintag{1.3}\ub{\stag{1.3} Question}:  If $|X|$ is sufficiently large and
${\frak D} \subseteq {\frak C}^{\text{full}}$ (is symmetric), when is
it true that maj-cl$({\frak D}) = {\frak C}^{\text{full}}$?
\bigskip

\definition{\stag{1.4} Definition}  1) Let Dis $= \text{ Dis}(X) = \{\mu:\mu$ a
distribution on ${\frak C}_X\}$; of course, ``$\mu$ a distribution on
${\frak C}$" means $\mu$ is a function from ${\frak C}$ into $[0,1]_{\Bbb
R}$ such that $\Sigma\{\mu(c):c \in {\frak C}\} = 1$.
\nl
2) For ${\Cal C} \subseteq {\frak C}$ and $\mu \in \text{Dis}({\frak C})$
let $\mu({\Cal C}) = \Sigma\{\mu(c):c \in {\Cal C}\}$ so 
$\mu({\Cal C}) \ge 0,\mu({\frak C}) = 1$. \nl
3) For ${\Cal D} \subseteq {\frak C}$ let Dis$_{\Cal D} = \{\mu \in
\text{ Dis}:\mu({\Cal D}) =1\}$. \nl
4) Let pr$({\frak C}) = \{\bar t:\bar t = \langle t_{x,y}:x \ne y
\in X \rangle$ such that $t_{x,y} \in [0,1]_{\Bbb R}$ and $t_{y,x} =
1-t_{x,y}\}$, we may write $\bar t(x,y)$ instead of $t_{x,y}$; pr 
stands for probability. \nl
5) For $T \subseteq \text{ pr}({\frak C})$ let pr-cl$(T)$ be the convex hull
of $T$. 
\nl
6) For $d \in {\frak C}$ let $\bar t[d] = \langle t_{x,y}[d]:x \ne y
\in X \rangle$ be defined by $[t_{x,y}[d] = 1 \Leftrightarrow d\{x,y\} = y
\Leftrightarrow t_{x,y}[d] \ne 0]$ when 
$\{x,y\} \in \text{ Dom}(d)$ and $t_{x,y}[d] = \frac 12 =
t_{y,x}[d]$ if $x \ne y \in X,\{x,y\} \notin \text{ Dom}(d)$.
\nl
7) Let pr-cl$({\Cal D})$ for ${\Cal D} \subseteq {\frak C}$ be
pr-cl$(\{\bar t[d]:d \in {\Cal D}\})$ and let prd$({\Cal D}) = \{\bar t[c]:c
\in {\Cal D}\}$. \nl
8) For ${\Cal C} \subseteq {\frak C}$ we let sym-$c \ell({\Cal C})$ be
the minimal ${\Cal D} \subseteq {\frak C}$ which is symmetric and
includes ${\Cal C}$.  For $T \subseteq \text{ pr}({\frak C})$ 
let maj$(T) = \{c \in {\frak C}$: for some $\bar t \in T$ we have
$c = \text{ maj}(\bar t)\}$, see below, and for 
${\Cal D} \subseteq {\frak C}$ let maj-cl$({\Cal D}) = 
\text{ maj(pr-cl}({\Cal D})$).
\nl
9) For $\bar t \in \text{ pr}({\frak C})$ we define maj$(\bar
t)$ as the $c \in {\frak C}$ such that $c\{x,y\} = y \Leftrightarrow
t_{x,y} > \frac 12$.
\enddefinition
\bigskip

\proclaim{\stag{1.8} Claim}  
1) For $d \in {\frak C}$ we have $\bar t[d] \in \,{\text{\rm pr\/}}({\frak
C})$. \nl
2) For ${\Cal D} \subseteq {\frak C}$ we have 
${\text{\rm Dis\/}}_{\Cal D} \subseteq { \text{\rm Dis\/}}$. \nl
3) ${\text{\rm prd\/}}({\frak C}) = { \text{\rm pr\/}}({\frak C})$ 
and if ${\Cal D} \subseteq {\frak C}$ \ub{then} 
${\text{\rm prd\/}}({\Cal D}) \subseteq {\text{\rm pr-cl\/}}({\Cal D})
\subseteq { \text{\rm Dis\/}}$. \nl
4) If ${\Cal C} \subseteq {\frak C}$ \ub{then} ${\Cal C} 
\subseteq { \text{\rm sym\/}}-c \ell({\Cal C}) \subseteq {\frak C}$. \nl
5) If $T \subseteq { \text{\rm pr\/}}({\frak C})$ \ub{then} 
${\text{\rm maj\/}}(T) \subseteq {\frak C}$. \nl
6) For ${\Cal D} \subseteq {\frak C}$ the two definitions of
${\text{\rm maj\/}}-c \ell({\Cal D}) = {\text{\rm maj(pr-cl\/}}({\Cal
D}))$ in \scite{1.4}(8) and \scite{1.2} are equivalent.
\endproclaim
\bigskip

\demo{Proof}  Obvious.

Kalai showed that not everything is possible.
\enddemo
\bigskip

\proclaim{\stag{1.9} Claim}  (G. Kalai) If ${\Cal C} \subseteq 
{\frak C}^{\text{full}}$ and 
for every $c \in {\Cal C}$ and $x \in X$, the in-valency and
out-valency are equal, (i.e., ${\text{\rm val\/}}_c(x) = (|X|-1)/2$,
see below) \ub{then} every $d \in {\frak C}^{\text{full}} \cap 
{\text{\rm maj-cl\/}}({\Cal C})$ satisfies:
\mr
\item "{$(*)$}"  if $\emptyset \ne Y \subsetneqq X$ then the
directed graph {\rm Tor}$(d)$ satisfies: there are
edges from $Y$ to $X \backslash Y$ and from $X \backslash Y$ to $Y$.
\endroster
\endproclaim
\bigskip

\demo{Proof}  See \scite{b.2}(1),(2) (and not used earlier).
\enddemo
\bigskip 

\definition{\stag{1.10} Definition}  1) For $d \in {\frak C}$ and $x
\in X$ let val$_d(x)$, the valency of $x$ for $d$ be
$|\{y:y \in X,y \ne x,d\{x,y\} = y\}| + |\{y:y \in X,y
\ne x,d\{x,y\}$ not defined$\}|/2$, so if $d$ is full the second term
disappears.  Let val$^+_d(x) = 
|\{y:y \in X,y \ne x$ and $d\{x,y\} = y\}|$ so val$^+_d(x) 
\in \{0,\dotsc,\bold n-1\}$. 
\nl
We also call val$^+_d(x)$ the out-valency \footnote{natural under the
tournament interpretation} of $x$ in $d$ and also call it
val$^{+1}_d(x)$ and we let $|\{y:y \in X,y\ne x$ and
$d\{y,x\} = x\}|$ be the in-valency of $x$ and denote it by
val$^{-1}_d(x)$; note that if $d$ is full (i.e., 
$\in {\frak C}^{\text{full}}_X$), then val$^{-1}_d(x) = \bold n -
\text{ val}^{+1}_d(x)-1$ and val$_d(x) = \text{ val}^+_d(x)$. \nl
2) For $d \in {\frak C}$ let Val$(d) = \{\text{val}_d(x):x \in X\}$.
\nl
3) For $d \in {\frak C}$ and $\ell \in \{0,1\}$ let $V_\ell(d) =
\{(\text{val}_d(x_0)$,val$_d(x_1)):x_0 \ne x_1 \in X$ and
$d\{x_0,x_1\} = x_\ell\}$. \nl
4) For $d \in {\frak C}$ and $\ell \in \{0,1\}$ let $V^*_\ell(d) =
\{\bar k-(\ell,1-\ell):\bar k \in V_\ell(d)\}$ and let $V^*(d) = V^*_0(d) \cup
V^*_1(d)$.
\nl
5) For $c \in {\frak C}$ let ${\text{\rm dual\/}} (c) \in {\frak C}$
have the same domain as $c$ and satisfy ${\text{\rm dual\/}}(c)\{x,y\} 
\in \{x,y\} \backslash \{c\{x,y\}\}$ when defined; 
similarly $\bar t' = \text{ dual}
(\bar t)$ for $\bar t \in \text{ pr}({\frak C})$ means that $t'_{x,y}
= 1 - t_{x,y}$. \nl
6) Let $V_{1/2}(d) = \{(\text{val}_d(x_0),\text{val}_d(x_1)):x_0 \ne
x_1 \in X$ and $d\{x_0,x_1\}$ is not defined$\}$ and $V^*_{1/2}(d) =
V_{1/2}(d)$. 
\enddefinition
\bigskip

\proclaim{\stag{1.11} Claim}  1) 
\mr
\item "{$(\alpha)$}"  $c_1 \in { \text {\rm sym-cl\/}}\{c_2\} 
\,{\text {\rm \ub{iff} dual\/}}(c_1) \in 
{ \text{\rm sym-cl\/}}\{\text{\rm dual}(c_2)\}$ 
\sn
\item "{$(\beta)$}"    $c_1 \in { \text{\rm maj-cl\/}}
({\text{\rm sym-cl\/}}\{c_2\}) 
\,{\text {\rm \ub{iff}\, dual\/}}(c_1) \in 
{\text{\rm maj-cl\/}}({\text{\rm sym-cl\/}}\{{\text{\rm dual\/}}(c_2)\})$.
\ermn
2) $(k_0,k_1) \in V_0(d) \Leftrightarrow (k_1,k_0) \in V_1(d)$ and
   $(k_0,k_1) \in V^*_0(d) \Leftrightarrow (b_1,b_0) \in V^*_1(d)$.
\nl
3) ``$c_1 \in { \text{\rm sym-cl\/}}\{c_2\}$" is an equivalence
relation on ${\frak C}$ and it implies $V_\ell(c_1) = V_\ell(c_2)$ for
$\ell=0,1$.  
\endproclaim
\bigskip

\demo{Proof}  Easy.
\enddemo
\newpage

\head {\S2 When every majority choice is possible: a characterization} \endhead  \resetall \sectno=2
 \spuriousreset
\bigskip

The following is the main part of the solution
(probably $(c) \Leftrightarrow (g)$ is the main conclusion here).
\proclaim{\stag{2.1} Main Claim}   Assume that  ${\frak D}
\subseteq {\frak C}^{\text{full}}$ which is symmetric and non-empty, 
(i.e., ${\frak D}$ is a non-empty set of choice functions on $\binom X2$ 
closed under permutation on $X$) and for simplicity assuming that 
${\frak D} = { \text{\rm sym-cl\/}}(d^*)$ for 
any $d^* \in {\frak D}$.  \ub{Then} the 
following conditions on ${\frak D}$ are 
equivalent, where $x,y$ vary on distinct members of $X$:
\mr
\widestnumber\item{$(b)_{x,y}$}
\item "{$(a)$}"  {\rm maj-cl}$({\frak D}) \supseteq {\frak
C}^{\text{full}}_X$
\sn
\item "{$(a)'$}"  ${\text{\rm maj-cl\/}}({\frak D}) = {\frak C}$
\sn
\item "{$(b)_{x,y}$}"  there is $\bar t \in { \text{\rm pr-cl\/}}({\frak D})
\subseteq {\text{\rm pr\/}}({\frak C})$ such that 
{\roster
\itemitem{ $(i)$ }  $t_{x,y} > \frac 12$
\sn
\itemitem{ $(ii)$ }  $\{x,y\} \ne \{u,v\} \in \binom X2 
\Rightarrow t_{u,v} = \frac 12$
\endroster}
\item "{$(c)$}"  for any $c \in {\frak C}^{\text{full}}$ we can 
find a finite set
$J$ and sequence $\langle d_j:j \in J \rangle$ such that $d_j \in
{\frak D}$ and: if $u \ne v \in X$ then 
\nl
$c\{u,v\} = v \Leftrightarrow |\{j \in J:d_j\{u,v\} = v\}| > |J|/2$
\sn
\item "{$(c)'$}"  like clause (c) for $c \in {\frak C}$
\sn
\item "{$(d)$}"  $(\frac 12,\frac 12)$ belongs to 
${\text{\rm Pr\/}}_{> \frac 12}({\frak D})$, see Definition
\scite{2.2} below
\sn
\item "{$(e)$}"  $(\frac 12,\frac 12) \in {\text{\rm Pr\/}}_{\ne 1/2}({\frak
D})$
\sn
\item "{$(f)$}"  $(\frac{\bold n}2 - 1,\frac{\bold n}2 -1)$ can be
represented as $r^*_0 \times \bar s_0 + r^*_1 \times \bar s_1$ where 
{\roster
\itemitem{ $(*)(i)$ }  $r^*_0,r^*_1 
\in [0,1]_{\Bbb R} \backslash \{\frac 12\}$
\sn
\itemitem{ $(ii)$ }  $1 = r^*_0 + r^*_1$
\sn
\itemitem{ $(iii)$ }  for $\ell =0,1$ the pair $\bar s_\ell \in \Bbb R
\times \Bbb R \text{ belongs to the convex hull of } V^*_\ell(d^*)$
\nl
$\text{for some } d^* \in {\frak D}$, see Definition \scite{1.10}(4), 
but recall that by a hypothesis of the claim, the
choice of $d^*$ is immaterial
\endroster}
\sn
\item "{$(g)$}"  for some $(d^* \in {\frak D}$ and) $x \in X$ we have
${\text{\rm val\/}}_{d^*}(x) \ne \frac{\bold n-1}2$. 
\endroster
\endproclaim
\bigskip

\demo{Proof}  \ub{$(b)_{x,y} \Leftrightarrow (b)_{x',y'}$}:

(So $x,y,x',y' \in X$ and $x \ne y,x' \ne y'$).  
Trivial as ${\frak D}$ is closed under permutations of $X$ hence so is
pr-cl$({\frak D})$.
\mn
\ub{$(b)_{x,y} \Rightarrow (a)'$}:

Let $c \in {\frak C}$. \nl
Let $\{(u_i,v_i):i < i(*)\}$ without repetitions list the 
pairs $(u,v)$ of distinct
members of $X$ such that $c\{u,v\} = v$; clearly $i(*) \le
\binom{|X|}{2}$ and $c \in {\frak C}^{\text{full}} \Rightarrow i(*) =
\binom {|X|}{2}$.   For
each $i < i(*)$ as $(b)_{x,y} \Rightarrow (b)_{u_i,v_i}$ clearly
there is $\bar t^i \in \text{ pr-cl}({\frak D})$
such that

$$
t^i_{u_i,v_i} > \frac 12 \text{ so } t^i_{v_i,u_i} = 1 - t_{u_i,v_i} <
\frac 12
$$

$$
\{u_i,v_i\} \ne \{u,v\} \in \binom X2 \Rightarrow t^i_{u,v} = \frac
12.
$$
\mn
Let $\bar t^* = \langle t^*_{u,v}:u \ne v \in X \rangle$ be defined
by

$$
t^*_{u,v} = \Sigma\{t^i_{u,v}:i < i(*)\}/i(*).
$$
\mn
As pr-cl$({\frak D})$ is convex and $i < i(*) \Rightarrow \bar t^i
\in \text{ pr-cl}({\frak D})$ clearly $\bar t^* \in \text{ pr-cl}
({\frak D})$.  Now for each $j < i(*),t^i_{u_j,v_j}$ is 
$\frac 12$ if $i \ne j$ and is $> \frac 12$ if $i=j$.   Hence
$t^*_{u_j,v_j}$ being the average of $\langle t^i_{u_j,v_j}:i < i(*)
\rangle$ is $> \frac 12$.  Hence $t^*_{v_j,u_j} = 1-t^*_{u_j,v_j} <
\frac 12$.  
So by the choice of $\langle(u_i,v_i):i < i(*) \rangle$ we have
\nl
$c\{u,v\} = v \Rightarrow t^*_{u,v} > \frac 12$ hence $c\{u,v\} = u
\Rightarrow t^*_{u,v} < \frac 12$.  Now lastly $c\{u,v\}$ undefined
$\Rightarrow \dsize \bigwedge_i \, t^i_{u,v} = \frac 12 \Rightarrow
t^*_{u,v} = \frac 12$.
So $\bar t^*$ witness $c \in$ maj-cl$({\frak D})$ as required in clause $(a)'$.
\mn
\ub{$(a)' \Rightarrow (a)$}: 

Trivial.
\mn
\ub{$(a) \Rightarrow (b)_{x,y}$}: 

By clause (a), for every $d \in {\frak C}^{\text{full}}$ 
there is $\langle r_c:c \in
{\frak D} \rangle$ as in Definition \scite{1.2}, hence for some
$\varepsilon_d > 0,u \ne v \in X \wedge d\{u,v\} = v
\Rightarrow \frac 12 + \varepsilon_d
< \Sigma\{r_c:c \in {\frak D}$ and $c\{u,v\} = v\}$.
Hence $\varepsilon = \text{ Min}\{\varepsilon_d:d \in {\frak D}\}$ 
is a real $>0$. 
\nl
Let $T = \{\bar t:\bar t \in \text{ pr-cl}({\frak D})$ 
and $t_{x,y} \ge \frac 12 + \varepsilon\}$, so
\mr
\item "{$(*)_1$}"  $T \ne \emptyset$
\nl
[Why?  By the choice of $\varepsilon$ and recall that ${\frak D}$ is symmetric]
\sn
\item "{$(*)_2$}"    $T$ is convex and closed.
\ermn
[Why?  Trivial.] \nl
For $\bar t \in T$ define 
\mr
\item "{$\boxtimes$}"  err$(\bar t) = \max\{|t_{u,v} - \frac 12|:u \ne
v \in X \text{ and } \{u,v\} \ne \{x,y\}\}$
\sn
\item "{$(*)_3$}"   if $\bar t \in T$, err$(\bar t) > 0$ \ub{then} we
can find $\bar t' \in T$ 
such that err$(\bar t') \le \text{ err}(\bar t)(1-\text{err}(\bar t))$ 
and $t'_{x,y} \ge (t_{x,y} + \frac 12 + \varepsilon)/2 \ge \frac 12 +
\varepsilon$.
\ermn
Why?  Choose $d \in {\frak C}^{\text{full}}$ such that $d\{x,y\} = y$ and

$$
u \ne v \in X \and \{u,v\} \ne \{x,y\} \and t_{u,v} > \frac 12
\Rightarrow d\{u,v\} = u
$$
\mn
(so if $t_{u,v} = t_{v,u} = \frac 12$ it does not matter what is
$d\{u,v\}$; such $d$ exists trivially).

So $d$ is ``a try to correct $\bar t$".

As we are assuming clause (a) and by the choice of $\varepsilon_d$,
we can find $\bar r^* = \langle r^*_c:c \in {\frak D} \rangle$ with 
$r^*_c \in [0,1]_{\Bbb R}$ and $1 =
\Sigma\{r^*_c:c \in {\frak D}\}$ such that  

$$
\frac 12 + \varepsilon_d < \Sigma\{r^*_c:c \in {\frak D} \text{ and } 
c\{x,y\} = y\}
$$
\mn
and if $u \ne v$ are from $X$ and $\{u,v\} \ne \{x,y\}$ then

$$
d\{u,v\} = v \Rightarrow \frac 12 < \Sigma\{r^*_c:c \in {\frak D}
\text{ and } c\{u,v\} = v\}
$$
hence

$$
d\{u,v\} = u \Rightarrow \frac 12 > \Sigma\{r^*_c:c \in {\frak D}
\text{ and } c\{u,v\} = v\}.
$$
\mn
By the choice of $\varepsilon$ \wilog \, 
$\frac 12 + \varepsilon < \Sigma\{r^*_c:c \in {\frak D}$ and $c\{x,y\}
= y\}$. \nl
Let $\bar s = \langle s_{u,v}:u \ne v \in X \rangle$ be defined by
$s_{u,v} = \Sigma\{r^*_c:c\{u,v\} = v\}$, so
\mr
\item "{$\circledast_1$}"   $(i) \quad \bar s \in \text{ pr-cl}({\frak D})$
\sn
\item "{${{}}$}"  $(ii) \quad s_{x,y} > \frac 12 + \varepsilon$ (so $s_{y,x} <
\frac 12$)
\sn
\item "{${{}}$}"  $(iii) \quad$ if 
$t_{u,v} > \frac 12$ and $u \ne v \in X,\{u,v\}
\ne \{x,y\}$ then $d\{u,v\} = u$ hence 
\nl

\hskip25pt $s_{u,v} < \frac 12$
\sn
\item "{${{}}$}"  $(iv) \quad$ if 
$t_{u,v} < \frac 12$ and $u \ne v \in X,\{u,v\}
\ne \{x,y\}$ then $d\{u,v\} = v$ hence
\nl

\hskip25pt  $s_{u,v} > \frac 12$.
\ermn
Choose $\delta \in (0,1)_{\Bbb R}$ as err$(\bar t)$.
Let $\bar t' = (1 - \delta) \bar t + \delta \bar s$, i.e. $t'_{u,v} =
((1-\delta)t_{u,v} + \delta s_{u,v})$ so clearly
\mr
\item "{$\circledast_2$}"  $(i) \quad \bar t' \in \text{ pr-cl}({\frak D})$
\sn
\item "{${{}}$}"  $(ii) \quad 
t'_{x,y} \ge \frac 12 + \varepsilon$ 
\sn
\item "{${{}}$}"  $(iii) \quad$ if $ u \ne v \in X,\{u,v\} \ne \{x,y\}$ then
$|t'_{u,v} - \frac 12| \le \text{ err}(\bar t)(1 - \text{ err}(\bar t))$
\sn
\item "{${{}}$}"  $(iv) \quad \bar t' \in T$.
\ermn
[Why?  Clause (i) as pr-cl$({\frak D})$ is convex.  Clause (ii) as
easily $t'_{x,y} = ((1-\delta)t_{x,y} + \delta s_{x,y})$, 
but $t_{x,y} \ge \frac 12 +
\varepsilon$ as $\bar t \in T$ and $s_{x,y} \ge \frac 12 +
\varepsilon$ by $\circledast_1(ii)$.  Now the main point, for clause (iii) note
that $t'_{u,v} - \frac 12 = -(t'_{v,u} - \frac 12)$ so as $d \in
{\frak C}^{\text{full}}$ \wilog \, $d\{u,v\} = u$ hence 
$t_{u,v} \ge \frac 12$, hence by the choice
of $d$ we have  $s_{u,v} \le \frac 12$ and both are in $[0,1]_{\Bbb R}$ and:

$$
\align
|t'_{u,v} - \frac 12| &= |(1-\delta)(t_{u,v} - \frac 12) +
\delta(s_{u,v} - \frac 12)| \le 
\underset {s \in [0,\frac 12]_{\Bbb R}} \to {\text{max}}|(1- \delta)
(t_{u,v} - \frac 12) + \delta(s - \frac 12)| \\
  &= \text{ Max}\{|(1-\delta)(t_{u,v} - \frac 12) + \delta(\frac 12 -
\frac 12)|,|(1-\delta)(t_{u,v} - \frac 12) + \delta(0 - \frac 12)|\} \\
  &=\text{ Max}\{|(1-\delta)(t_{u,v} - \frac 12),|(1- \delta)(t_{u,v}
- \frac 12) - \frac 12 \delta|\}  \\
  &\le \text{ Max}\{(1-\delta)(t_{u,v} - \frac 12),(1-\delta)(t_{u,v} 
- \frac 12),\frac 12 \delta\} \\
  &\le \text{ Max}\{(1-\delta)\text{ err}(\bar t),(1-\delta)
\text{ err}(\bar t), 
\frac 12 \delta\} \le \text{ err}(\bar t)(1 - \text{ err}(\bar t))
\endalign
$$
\mn
(recalling $\delta = \text{ err}(\bar t) \in [0,\frac 12]_{\Bbb R}$ so
$\frac 12 \delta \le \text{ err}(\bar t)(1 - \text{ err}(\bar t))$ 
as required), so clause (iii) holds.

Clause (iv) follows.  So $\circledast_2$ holds.]

So we are done proving $(*)_3$.

As $T$ is closed (and is included in a $\{\bar t:\bar t = \langle
t_{u,v}:u \ne v \in X \rangle$ and $0 \le t_{u,v} \le 1\}$ which is
compact), clearly there
is $\bar t \in T$ such that $u \ne v \in X \and \{u,v\} \ne \{x,y\}
\Rightarrow t_{u,v} = \frac 12$ as required in part (ii) of $(b)_{x,y}$.
\mn
\ub{$(c)' \Rightarrow (c)$}:

Trivial.
\mn
\ub{$(c) \Rightarrow (a)$}:

Let $d^* \in {\frak C}^{\text{full}}$ and 
let $\langle c_j:j \in J \rangle$ witness clause (c) for $d^*$. 
\nl
Let $r_c = |\{j \in J:c_j = c\}|/|J|$ now $\langle r_c:c \in {\frak
D} \rangle$ witness clause (a), i.e., witness that $d^* \in$  
maj-cl$({\frak D})$.
\mn
\ub{$(a) \Rightarrow (c)$}:

Let $d^* \in {\frak C}^{\text{full}}$ and let $\langle r_c:c \in {\frak D}
\rangle$ be as guaranteed for $d^*$ by clause $(a)$.  Let $n(*) > 0$ be large
enough such that $|\frac{1}{2} - r_c| > \frac{1}{n(*)}$ for $c \in
{\frak K}$ 
and for $c \in {\frak D}$ let $k_c \in \{0,\dotsc,n(*)-1\}$ be
such that $c \in {\frak D} \Rightarrow k_c \le n(*) \times r_c < k_c +
1$; note that $k_c$ exists as $r_c \in [0,1]_{\Bbb R}$.
As $\dsize \sum_c \frac{k_c}{n(*)} \le 1 \le \dsize \sum_c
\frac{k_c+1}{n(*)}$, we can choose $m_c \in \{k_c,k_c+1\}$ such that $r'_c =
\frac{m_c}{n(*)}$ satisfies $\Sigma\{r'_c:c \in {\frak D}\} = 1$.
Let $J = \{(c,m):c \in {\frak D}$ and $m \in \{1,\dotsc,m_c\}\}$ and
we let $c_{(d,m)} = d$ for $(d,m) \in J$.  
Now the ``majority" of $\langle c_t:t \in J \rangle$, 
see Definition \scite{1.2}, choose $d^*$ so clause (c)
holds.
\mn
\ub{$(a)' \Rightarrow (c)'$}:

Similar, using: if a finite set of equalities and inequalities with
rational coefficients is solvable in $\Bbb R$ then it is solvable in 
$\Bbb Q$.
\enddemo
\bn
Before we deal with clauses (d),(e),(f) and (g) of \scite{2.1}, we define
\definition{\stag{2.2} Definition}  1) For ${\Cal D} \subseteq 
{\frak C}^{\text{full}}$ and $A \subseteq [0,1]_{\Bbb R}$ let 
Pr$_A({\Cal D})$ be the set of pairs $(s_0,s_1)$ of real
numbers $\in [0,1]_{\Bbb R}$ such that 
for some $\bar t \in \text{ pr-cl}({\Cal D})$ and $x \ne y \in X$ and
$a \in A$ we have $\bar t = \bar t \langle x,y,a,s_0,s_1 \rangle$
where \nl
2) $\bar t = \bar t\langle x,y,a,s_0,s_1 \rangle$ where
$x \ne y \in X,a \in [0,1]_{\Bbb R}$ and $s_0,s_1 \in [0,1]_{\Bbb R}$ and
$\bar t = \langle t_{u,v}:u \ne v \in X \rangle \in 
\text{ pr}({\frak C})$ is defined by 
\mr
\item "{$(\alpha)$}"  $t_{x,y} = a$
\sn
\item "{$(\beta)$}"  if $z \in X \backslash \{x,y\}$ \ub{then}
$t_{y,z} = s_1$ (hence $t_{z,y} = 1 - s_1$)
\sn
\item "{$(\gamma)$}"  if $z \in X \backslash \{x,y\}$ \ub{then}
$t_{x,z} = s_0$ (hence $t_{z,x} = 1-s_0$)
\sn
\item "{$(\delta)$}"  if $z_1 \ne z_2 \in X \backslash \{x,y\}$
\ub{then} $t_{z_1,z_2} = \frac 12$.
\ermn
3) In Pr$_A({\Cal D})$ we may replace $A$ by $1,0,\ne\frac 12,> \frac
12,< \frac 12$ if $A$ is
$\{1\},\{0\},[0,1]_{\Bbb R} \backslash \{\frac 12\},(\frac 12,1]_{\Bbb
R},[0,\frac 12)_{\Bbb R}$ respectively. 
\nl
4) For $\ell \in \{0,1\}$ and 
${\Cal D} \subseteq {\frak C}^{\text{full}}$ 
let Prd$_\ell({\Cal D})$ be the set
of pairs $\bar s = (s_0,s_1)$ of real (actually rational) 
numbers $\in [0,1]_{\Bbb R}$
such that for some $c \in {\Cal D}$ and $x \ne y \in X$ we have $\bar
s = \bar s^{c,x,y} = (s^{c,x,y}_0,s^{c,x,y}_1)$ where
\mr
\widestnumber\item{$(iii)$}
\item "{$(i)$}"  $s^{c,x,y}_1 = 
|\{z:z \in X \backslash \{x,y\} \text{ and } c\{y,z\} = z\}|/(\bold n - 2)$
\sn
\item "{$(ii)$}"  $s^{c,x,y}_0 = 
|\{z:z \in X \backslash \{x,y\} \text{ and } c\{x,z\} = z\}|/(\bold n - 2)$
\sn
\item "{$(iii)$}"  $\ell = 1 
\Leftrightarrow \ell \ne 0 \Leftrightarrow c\{x,y\} = y$.
\ermn
5) For ${\Cal D} \subseteq {\frak C}^{\text{full}}$ let
Prd$({\Cal D})$ be Prd$_0({\Cal D}) \cup \text{ Prd}_1({\Cal D})$.
\enddefinition
\bigskip

\proclaim{\stag{2.3} Claim}   Let ${\Cal D} \subseteq {\frak C}^{\text{full}}$
\mr
\item  ${\text{\rm Pr\/}}_{A_1}
({\Cal D}_1) \subseteq { \text{\rm Pr\/}}_{A_2}({\Cal D}_2)$
if $A_1 \subseteq A_2 \subseteq [0,1]_{\Bbb R}$ and ${\Cal D}_1
\subseteq {\Cal D}_2 \subseteq {\frak C}^{\text{full}}$
\sn
\item   ${\text{\rm Pr\/}}_A({\Cal D})$ is a convex subset of
$[0,1]_{\Bbb R} \times [0,1]_{\Bbb R}$ when $A$ is a convex subset of
$[0,1]_{\Bbb R}$
\sn
\item  ${\text{\rm Prd\/}}_\ell
({\Cal D})$ is finite and its convex hull is $\subseteq
{\text{\rm Pr\/}}_\ell({\Cal D})$, increasing with ${\Cal D}(\subseteq
{\frak C}^{\text{full}})$ for $\ell =0,1$
\sn
\item  For $x \ne y \in X$ and $c \in {\frak C}^{\text{full}}$ 
satisfying $\ell = 1
\Rightarrow c\{x,y\} = y$ and $\ell=0 \Rightarrow c\{x,y\} = x$
we have (see Definition \scite{1.10}(1))
$$
s^{c,x,y}_0 = ({\text{\rm val\/}}_c(x) - \ell)/(\bold n - 2)
$$

$$
s^{c,x,y}_1 = ({\text{\rm val\/}}_c(y) - (1 - \ell))/(\bold n - 2).
$$
\sn
\item  ${\text{\rm Pr\/}}_{A_1 \cup A_2}
({\Cal D}) = { \text{\rm Pr\/}}_{A_1}({\Cal D}) \cup
{ \text{\rm Pr\/}}_{A_2}({\Cal D})$, in fact ${\text{\rm Pr\/}}_A({\Cal D}) =
\cup\{{\text{\rm Pr\/}}_{\{a\}}({\Cal D}):a \in A\}$.
\endroster
\endproclaim
\bigskip

\demo{Proof}  Immediate; part (3) holds by \scite{2.3A}(1) below
concerning part (4) recall Definition \scite{1.10}(1).
\hfill$\square_{\scite{2.3}}$ 
\enddemo
\bigskip

\proclaim{\stag{2.3A} Claim}  1) If $x \ne y \in X$ and $c \in {\frak
C}$ and $\ell \in \{0,1\}$ satisfies $\ell = 1 \Rightarrow c\{x,y\}
= y$ and $\ell=0 \Rightarrow c\{x,y\} = x$ recalling Definition
\scite{1.4}(6) and \scite{1.1A}(1) we have: 
$\bar t \langle x,y,\ell,s^{c,x,y}_0,s^{c,x,y}_1
\rangle = \frac{1}{|\Pi_{x,y}|} \Sigma\{\bar t[\hat \pi(c)]:\pi \in
\Pi_{x,y}\}$ where $\Pi_{x,y} := \{\pi \in { \text{\rm Per\/}}(X):\pi(x) = x$
and $\pi(y) = y\}$ hence $|\Pi_{x,y}| = (\bold n -2)!$. 
\nl
2) If ${\Cal D} \subseteq {\frak C}$ is symmetric, $x \ne y \in X$ 
and $\bar t \in 
{\text{\rm pr-cl\/}}({\Cal D})$ and $\bar t^* = \Sigma\{\bar t^\pi:\pi
\in \Pi_{x,y}\}/|\Pi_{x,y}|$ where $\bar t^\pi = \langle t^\pi_{u,v}:u
\ne v \in X \rangle,t^\pi_{u,v} = t_{\pi^{-1}(u),\pi^{-1}(v)}$
\ub{then} $\bar t^* \in { \text{\rm pr-cl\/}}({\Cal D})$ and $(s_0,s_1) \in
{ \text{\rm Pr\/}}_{\{a\}}({\Cal D})$ where $a = t_{x,y},s_0 =
\Sigma\{t_{x,z}:z \in X \backslash \{x,y\}\}/(\bold n-2)$ and $s_1 =
\Sigma\{t_{y,z}:z \in X \backslash \{x,y\}\}/(\bold n-2)$.
\endproclaim
\bigskip

\demo{Proof}  Easy (in part (2), $\bar t^*$ witness $(s_0,s_1) \in
\text{ Pr}_{\{a\}}({\Cal D}))$.   \hfill$\square_{\scite{2.3A}}$
\enddemo
\bigskip

\proclaim{\stag{2.4} Claim}  For any symmetric non-empty ${\Cal D} \subseteq
{\frak C}^{\text{full}}$ (i.e., closed under permutations of $X$): \nl
1) For $\ell \in \{0,1\}$, the set ${\text{\rm Pr\/}}_\ell
({\Cal D})$ is the convex
hull of ${\text{\rm Prd\/}}_\ell({\Cal D})$ in $\Bbb R \times \Bbb R$. 
\nl
2) ${\text{\rm Pr\/}}_{[0,1]}({\Cal D})$ is the convex hull of 
${\text{\rm Prd\/}}({\Cal D})$. 
\nl
3) Let $a \in [0,1]_{\Bbb R}$ and 
$s^*_0,s^*_1 \in [0,1]_{\Bbb R}$.  \ub{Then} $\bar t^* = \bar
t\langle x,y,a,s^*_0,s^*_1 \rangle \in {\text{\rm pr-cl\/}}({\Cal
D})$ \ub{iff} we can find $\langle r_{\bar s,\ell}:\ell \in \{0,1\}$
and $\bar s \in {\text{\rm Prd\/}}_\ell({\Cal D})\rangle$ such that
$r_{\bar s,\ell} \in [0,1]_{\Bbb R}$ and $1 = \Sigma\{r_{\bar
s,\ell}:\ell \in \{0,1\},\bar s \in 
{\text{\rm Prd\/}}_\ell({\Cal D})\}$ and $(s^*_0,s^*_1) = \Sigma\{r_{\bar
s,\ell} \times \bar s:\ell \in \{0,1\},\bar s \in { \text{\rm
Prd\/}}_\ell({\Cal D})\}$ and $a = \Sigma\{r_{\bar s,1}:
\bar s \in {\text{\rm Prd\/}}_1({\Cal D})\}$.
\endproclaim
\bigskip

\demo{Proof}  1) By \scite{2.3}(3) we have one inclusion.

For the other direction assume $(s^*_0,s^*_1) \in \text{
Pr}_\ell({\Cal D})$ and we should prove that the pair $(s^*_0,s^*_1)$ 
belongs to the convex hull of Prd$_\ell({\Cal D})$.
Fix $x \ne y \in X$ and let $\bar t^* = \bar t \langle
x,y,\ell,s^*_0,s^*_1 \rangle$ so
\mr
\item "{$\boxtimes_1$}"  $\bar t^* = 
\langle t^*_{u,v}:u \ne v \in X \rangle$ is defined as
follows $t^*_{x,y} = \ell,t^*_{u,v} = \frac 12$ if 
$u \ne v \in X \backslash \{x,y\},t^*_{y,z} =
s^*_1$ if $z \in X \backslash \{x,y\},t^*_{x,z} = s^*_0$ if $z \in X
\backslash \{x,y\}$.
\ermn
As $(s^*_0,s^*_1) \in \text{ Pr}_\ell({\Cal D})$ by 
Definition \scite{2.2} we know that $\bar t^* \in \text{ pr-cl}({\Cal D})$ and
let $\bar r = \langle r_c:c \in {\Cal D} \rangle$ be such that 
\mr
\item "{$\boxtimes_2$}"  $x,y,\bar r$ witness that 
$\bar t^* \in \text{ pr-cl}({\Cal D})$, so
$r_c \ge 0$ and $1 = \Sigma\{r_c:c \in {\Cal D}\}$ and

$$
\bar t^* = \Sigma\{r_c \times \bar t[c]:c \in {\Cal D}\}.
$$
\ermn
As $t^*_{x,y} = \ell$, necessarily
\mr
\item "{$\boxtimes_3$}"   $r_c \ne 0 \Rightarrow c \in 
{\Cal D}^\ell_{x,y} := \{c \in {\Cal D}:(\ell =1 \Rightarrow 
c\{x,y\} = y)$ and $(\ell=0 \Rightarrow c\{x,y\} = x)\}$.
\ermn
To make the rest of the proof also a proof of part (3) let $a = \ell$
(as the real number $a$ may be $\ne 0,1$; in any case we use
$m \in \{0,1\}$ below).
\nl
Let $\Pi_{x,y} = \{\pi \in \text{ Per}(X):\pi(x) = x,\pi(y) = y\}$ and
recall that for $\pi \in \text{ Per}(X),\hat \pi$ is the 
permutation of ${\frak C}$ which $\pi$ induces, defined in
\scite{1.1A}, so $\hat \pi$ maps ${\Cal D}^\ell_{x,y}$ 
onto ${\Cal D}^\ell_{x,y}$ if $\pi \in \Pi_{x,y}$ recalling we have
assumed that ${\Cal D}$ is symmetric.  Clearly
$|\Pi_{x,y}| = (\bold n - 2)!$; recall that $\bar t^{**} = \bar
t^*[\hat \pi[c]]$ if $t^{**}_{u,v} = t^*_{\pi(u),\pi(v)}$.
\nl
For $(s_0,s_1) \in \text{ Prd}_\ell({\Cal D}) \subseteq 
\{(\frac{m_1}{(\bold n-2)!},\frac{m_2}{(\bold n-2)!}):m_1,m_2 \in 
\{0,1,\dotsc,(\bold n-2)!\}\}$ let
\mr
\item "{$\boxtimes_4$}"   $r^*_{(s_0,s_1)} =
\Sigma\{r_c:\bar s^{c,x,y} = (s_0,s_1)\}$ and for $m \in \{0,1\}$ we 
let $r^*_{(s_0,s_1),m} = \Sigma\{r_c:\bar s^{c,x,y} = (s_0,s_1)$ and
$m=1 \Rightarrow c\{x,y\} = y$ and $m=0 \Rightarrow c\{x,y\} = x\}$.
\ermn
Clearly $\pi \in \Pi_{x,y} \Rightarrow \langle t^*_{\pi(u),\pi(v)}:u \ne
v \in X \rangle = \bar t^*$, just check the definition, hence
(by the beginning of this sentence; by the equation in $\boxtimes_2$;
by arithmetic; by \scite{2.3A}(1); as $\langle \{c \in {\Cal D}:m=1
\Rightarrow c \{x,y\} = y$ and $m=0 \Rightarrow c\{x,y\} = x$ and
$\bar s^{c,x,y} = (s_0,s_1)\}:m \in\{0,1\}$ and $(s_0,s_1) \in \text{
Prd}_m({\Cal D})\rangle$ is a partition of ${\Cal D}$; by the choice
of $r^*_{(s_0,s_1)}$ in $\boxtimes_4$) we have:
\mr
\item "{$\boxtimes_5$}"  $\bar t^* = \dsize \frac{1}{| \Pi_{x,y}|} \dsize
\sum_{\pi \in \Pi_{x,y}} \langle t^*_{\pi(u),\pi(v)}:u \ne v \in X
\rangle = \frac{1}{|\Pi_{x,y}|} \dsize \sum_{\pi \in \Pi_{x,y}} \,
\dsize \sum_{c \in {\Cal D}} r_c \times \bar t^*[\hat \pi(c)] =$ \nl 
$\dsize \sum_{c \in {\Cal D}} \, r_c(\frac{1}{|\Pi_{x,y}|} \, \dsize
\sum_{\pi \in \Pi_{x,y}} \bar t^*[\hat \pi(c)]) =
\dsize \sum_{c \in {\Cal D}} r_c \times \bar t \langle
x,y,a,s^{c,x,y}_0,s^{c,x,y}_1 \rangle =$ \nl
$\dsize \sum_{m \in \{0,1\}} \,\, 
\dsize \sum_{(s_0,s_1) \in \text{ Prd}_m({\Cal D})} 
\, (\Sigma\{r_c:c \in {\Cal D},m=1 \Rightarrow c\{x,y\} = y$ \nl

\hskip110pt and $m=0 \Rightarrow c\{x,y\} = x$ \nl

\hskip110pt and $\bar s^{c,x,y} = (s_0,s_1)\}) \times
\bar t \langle x,y,m,s_0,s_1) =$ \nl
$\dsize \sum_{m \in \{0,1\}} \, 
\dsize \sum_{(s_0,s_1) \in \text{ Prd}_m({\Cal D})}
r^*_{(s_0,s_1),m} \bar t\langle x,y,m,s_0,s_1 \rangle$.
\ermn
Now concentrate again on the case $a = \ell \in \{0,1\}$, so
$r^*_{\bar s,1-\ell} = 0$ by $\boxtimes_3$ and $r^*_{\bar s,\ell} = 
r^*_{\bar s}$.
So clearly
\mr
\item "{$\circledast_1$}"  $r^*_{\bar s} \ge 0$ \nl
[Why?  As the sum of non-negative reals]
\sn
\item "{$\circledast_2$}"   $1 =
\Sigma\{r^*_{\bar s}:\bar s \in \text{ Prd}_\ell({\Cal D})\}$ \nl
[Why?  As by Definition \scite{2.2}(2), $c \in {\Cal D} \and r_c > 0
\Rightarrow c \in {\Cal D}^\ell_{x,y} \Rightarrow \bar s^{c,x,y} \in
\text{Prd}_\ell({\Cal D})$ and the definition of $r^*_{\bar s}$]
\sn
\item "{$\circledast_3$}"  we have
{\roster
\itemitem{ $(\alpha)$ }  $z \in X \backslash \{x,y\} \Rightarrow s^*_1
= t^*_{y,z} = \Sigma\{r^*_{(s_0,s_1)} \times s_1:(s_0,s_1) \in \text{
Prd}_\ell({\Cal D})\}$,
\sn
\itemitem{ $(\beta)$ }  $z \in X \backslash \{x,y\} \Rightarrow s^*_0
= t^*_{x,z} = \Sigma\{r^*_{(s_0,s_1)} \times s_0:(s_0,s_1) \in \text{ Prd}_\ell
({\Cal D})\}$.
\nl
[Why?  By $\boxtimes_5$ and the definition of $\bar t\langle
x,y,m,s_0,s_1\rangle$.]
\endroster}
\ermn
So $\langle r^*_{\bar s}:\bar s \in \text{ Prd}_\ell({\Cal D}) \rangle$
witness that $(s^*_0,s^*_1) \in$ convex hull of Prd$_\ell({\Cal
D})$. \nl
2) Similar proof (and not used). 
\nl
3) One direction is as in \scite{2.3}(3).  For the other, by the
hypothesis, $\boxtimes_2$ in the proof of part (1) with $\ell$
replaced by $a$ holds.  So by the part of the proof of part (1) from
$\boxtimes_3$ till (and including) $\boxtimes_5$ we know that
$r^*_{\bar s,m}$ are defined and $\boxtimes_5$ holds.  So
\mr
\item "{$\boxdot_1$}"  $r^*_{\bar s,m} \ge 0$
\sn
\item "{$\boxdot_2$}"  $1= \Sigma\{r^*_{\bar s,m}:m \in \{0,1\}$ and
$\bar s \in$ Prd$_m({\Cal D})\}$
\sn
\item "{$\boxdot_3$}"  $(s^*_0,s^*_1) = \Sigma\{r^*_{\bar s,m} \times
\bar s:m \in \{0,1\},\bar s \in$ Prd$_m({\Cal D})\}$ \nl
[Why?  By $\boxtimes_5$.]
\sn
\item "{$\boxdot_4$}"  $a = \Sigma\{r^*_{\bar s,1}:\bar s \in$
Prd$_1({\Cal D})\}$.
\ermn
[Why?  Use $\boxtimes_5$, noting that $t^*_{x,y} = a$.]
\nl
So we are done.   \hfill$\square_{\scite{2.4}}$
\enddemo
\bn
\ub{Continuation of the proof of \scite{2.1}}:
\sn
\ub{$(d) \Leftrightarrow (b)_{x,y}$}:

Read Definition \scite{2.2}(1) 
(and the symmetry).  \hfill$\square_{\scite{2.1}}$
\bn
\ub{$(d) \Rightarrow (e)$}:

By \scite{2.3}(1).
\bn
\ub{$(e) \Rightarrow (d)$}:

Why?  If clause (e) holds, for some $a \in [0,1]_{\Bbb R} \backslash
\{\frac 12\}$ and $x \ne y \in X$ we have $\bar t^* =: \bar t \langle
x,y,a,\frac 12,\frac 12 \rangle \in$  pr-cl$({\frak D})$.  
If $a > \frac 12$ this witness $(\frac 12,\frac 12) \in \text{ Pr}_{>
1/2}({\Cal D})$, so assume $a < \frac 12$.  But trivially $\bar t
\langle y,x,1-a,\frac 12,\frac 12 \rangle$ is equal to $\bar t^*$ hence (as
in \scite{1.11}) is in pr-cl$({\frak D})$ and by symmetry we are done.
\bn
\ub{$(e) \Leftrightarrow (f)$}:

Clearly (e) means that
\mr
\item "{$(*)_0$}"  there are $r_c \in [0,1]_{\Bbb R}$ for 
$c \in {\frak D}$ such that 
$1 = \Sigma\{r_c:c \in {\frak D}\}$ and
$a \in [0,1]_{\Bbb R} \backslash \{\frac
12\}$ such that $\bar t \langle x,y,a, \frac 12,\frac 12 \rangle = 
\Sigma\{r_c \times \bar t[c]:c \in {\frak D}\}$. 
\ermn
By \scite{2.4}(3) we know that $(*)_0$ is equivalent to
\mr
\item "{$(*)_1$}"  there are $r_{\bar s,\ell} \in [0,1]_{\Bbb R}$ for $\bar
s \in \text{ Prd}_\ell({\frak D}),\ell \in \{0,1\}$ such that 
{\roster
\itemitem{ $(i)$ }  $1 = \Sigma\{r_{\bar s,\ell}:\bar s \in \text{
Prd}_\ell({\frak D}),\ell \in \{0,1\}\}$
\sn
\itemitem{ $(ii)$ }  $(\frac 12,\frac 12) = 
\Sigma\{r_{\bar s,\ell} \times \bar s:\bar s \in
\text{ Prd}_\ell({\frak D}),\ell \in \{0,1\}\}$
\sn
\itemitem{ $(iii)$ }  $\frac 12 \ne a = \Sigma\{r_{\bar s,1}:
\bar s \in \text{ Prd}_1({\frak D})\}$.
\endroster}
\ermn
But by \scite{2.3}(4) and Definition \scite{2.2}(4) for $\ell \in \{0,1\}$:

$$
\align
\text{Prd}_\ell({\frak D}) = \bigl\{&\bigl(\frac{\text{val}_c(x) - 
\ell}{\bold n -2}, \, \frac{\text{val}_c(y)-(1-\ell)}{\bold
n-2}\bigr):c \in {\frak D} \text{ and} \\
  &x \ne y \text{ and } (\ell = 1 \Rightarrow 
c\{x,y\} =y \text{ and } \ell=0 \Rightarrow c\{x,y\} = x)\bigr\}.
\endalign
$$
\mn
Let $d^* \in {\frak D}$, recall that ${\frak D} = 
\text{ sym-cl}(\{d^*\})$ by a hypothesis of \scite{2.1} and recall
$V_\ell(d^*) = \{(k_1,k_2):\text{ for
some } x_1 \ne x_2 \in X,k_1 = \text{ val}_{d^*}(x_1),k_2 = 
\text{ val}_{d^*}(x_2)$ and $d^*\{x_1,x_2\} = x_{\ell +1}\}$ for
$\ell=0,1$.  So $(*)_1$ means
(recalling the definition of Prd$_\ell({\frak D}))$
\mr
\item "{$(*)_2$}"  there is a sequence $\langle r_{\bar k,\ell}:\bar k
\in V_\ell(d^*)$ and $\ell \in \{0,1\} \rangle$ such that
{\roster
\itemitem{ $(i)$ }  $r_{\bar k,\ell} \in [0,1]_{\Bbb R}$ and
\sn
\itemitem{ $(ii)$ }  $1 = \Sigma\{r_{\bar k,\ell}:\bar k \in
V_\ell(d^*)$ and $\ell \in \{0,1\}\}$ and
\sn
\itemitem{ $(iii)$ }  $(\frac 12,\frac 12) = \Sigma \bigl\{
r_{(k_1,k_2),\ell} \times 
\bigl( \frac{k_1 - \ell}{\bold n-2}, \, 
\frac{k_2 -(1-\ell)}{\bold n-2} \bigr):\ell \in
\{0,1\}$ and \nl

\hskip65pt $(k_1,k_2) \in V_\ell(d^*)\}$
\sn
\itemitem{ $(iv)$ }  $\frac 12 \ne \Sigma\{r_{\bar k,1}:\bar k \in
V_1(d^*)\}$.
\endroster}
\ermn
Let us analyze $(*)_2$.  
Let $r^*_\ell = \Sigma\{r_{\bar k,\ell}:\bar k \in V_\ell(d^*)\}$ for
$\ell \in \{0,1\}$.  So
\mr
\item "{$\circledast_1$}"  $r^*_\ell \in 
[0,1]_{\Bbb R}$ and $1 = r^*_0 + r^*_1$.
\ermn
Now clause $(iii)$ of $(*)_2$  means $(iii)_1 + (iii)_2$ where
\mr
\item "{$(iii)_1$}"  $\frac 12 = \frac{1}{\bold n-2} \bigl(
\Sigma\{r_{\bar k,\ell} \times k_1:\bar k \in V_\ell(d^*)$ and $\ell \in
\{0,1\}\} \bigr) -$ \nl

\hskip30pt $\frac{1}{\bold n-2} \, \Sigma\{r_{\bar k,\ell} \times \ell:
\bar k \in V_\ell(d^*),\ell \in \{0,1\}\}) = $ \nl

\hskip30pt $\frac{1}{\bold n -2} \, \Sigma\{r_{\bar k,\ell} \times k_1:
\bar k \in V_\ell(d^*)$ and $\ell \in \{0,1\}\} - \frac{r^*_1}{\bold n-2}$, \nl
i.e., 
\sn
\item "{$(iii)'_1$}"  $\frac{\bold n}{2} - 
(1 - r^*_1) = \frac{\bold n -2}{2} + r^*_1 = 
\Sigma\{r_{\bar k,\ell} \times k_1:\bar k \in V_\ell(d^*)$ and $\ell \in
\{0,1\}\}$
\nl
and
\sn
\item "{$(iii)_2$}"  $\frac 12 = \frac{1}{\bold n -2} \,
\Sigma\{r_{\bar k,\ell} \times k_2:\bar k \in V_\ell(d^*)$ and $\ell \in
\{0,1\}\} - \frac{1}{\bold n -2} \, \Sigma\{r_{\bar k,\ell}
\times(1-\ell):$ \nl

\hskip30pt $\bar k \in V_\ell(d^*)$ and $\ell \in \{0,1\}\} =$ \nl

\hskip30pt $\frac{1}{\bold n -2} \, 
\Sigma\{r_{\bar k,\ell} \times k_2:\bar k \in
V_\ell(d^*),\ell \in \{0,1\}\} - \frac{r^*_0}{\bold n -2}$, \nl
i.e., 
\sn
\item "{$(iii)'_2$}"   $\frac{\bold n}{2} -(1 - r^*_0) = 
\frac{\bold n -2}{2} + r^*_0 = 
\Sigma\{r_{\bar k,\ell} \times k_2:\bar k \in V_\ell(d^*)$ and $\ell \in
\{0,1\}\}$.
\ermn
Together $(iii)$ of $(*)_2$ is equivalent to 
\mr
\item "{$(iii)^+$}"  $\bigl(\frac{\bold n}{2} - (1-r^*_1),\frac{\bold
n}{2} - (1-r^*_0) \bigr) =
\Sigma\{r_{\bar k,\ell} \times \bar k:\bar k \in V_\ell(d^*)$ and $\ell \in
\{0,1\}\}$.
\ermn
Let $\bar s_\ell = \Sigma\{r_{\bar k,\ell} \times \bar k:
\bar k \in V_\ell(d^*)\} / r^*_\ell$ if $r^*_\ell > 0$ and any member
of conv$(V_\ell(d^*))$ if $r^*_\ell = 0$, so $(*)_2$ is equivalent to
($V_\ell(d^*)$ is from Definition \scite{1.10})
\mr
\item "{$(*)_3$}"  there are $\bar s_0,\bar s_1,r^*_0,r^*_1$ such that
{\roster
\itemitem{ $(i)$ }  $\bar s_\ell \in \text{ conv}(V_\ell(d^*))$ for
$\ell = 0,1$
\sn
\itemitem{ $(ii)$ }  $r^*_0,r^*_1 \in [0,1]_{\Bbb R}$ and $1 = r^*_0 +
r^*_1$
\sn
\itemitem{ $(iii)$ }  $(\frac{\bold n}{2} - (1-r^*_1),\frac{\bold n}{2} -
(1-r^*_0))$ is $r^*_0 \times \bar s_0 + r^*_1 \times \bar s_1$
\sn
\itemitem{ $(iv)$ }  $r^*_\ell \ne \frac 12$ (by clause $(iv)$ in $(*)_2$
above).
\endroster}
\ermn
Clearly $(*)_3(iii)$ is equivalent to
\mr
\item "{$(iii)'$}"  $(\frac{\bold n}{2} - 1,\frac{\bold n}{2} -
1)$ is $r^*_0 \times (\bar s_0 - (0,1)) + r^*_1 \times 
(\bar s_1 - (1,0)\bigr)$.
\ermn
So $(*)_3$ is equivalent to
\mr
\item "{$(*)_4$}"   clauses $(i),(ii),(iv)$ of $(*)_3$
and $(iii)'$ above holds.
\ermn
But recalling Definition \scite{1.10}(4) of $V^*_\ell(d^*)$, 
this is clause (f), so we are done proving $(e) \Leftrightarrow (f)$.

\sn
\ub{$(g) \Rightarrow (f)$}:  By \scite{2.5}, \scite{2.6}, \scite{2.7}
below (i.e., they show $(g) + \neg(f)$ lead to a contradiction).
\bn
\ub{$(f) \Rightarrow (g)$}:  It suffices to prove $\neg(g) 
\Rightarrow \neg(f)$.  This holds trivially as 
$\neg(g)$ implies $(s_0,s_1) \in \text{ conv}(V_\ell(d^*))
\Rightarrow s_0 = s_1$.

We have proved $(b)_{x,y} \Leftrightarrow (b)_{x',y'},(b)_{x,y} \Rightarrow
(a)' \Rightarrow (a) \Rightarrow (b)_{x,y},(c)' \Rightarrow
(c) \Rightarrow (a) \Rightarrow
(c),(a)' \Rightarrow (c)' \Rightarrow (c),(d) \Leftrightarrow 
(b)_{x,y},(d) \Rightarrow (e) \Rightarrow (d),(e)
\Leftrightarrow (f),(g) \Rightarrow (f) \Rightarrow (g)$, 
so we are done proving \scite{2.1}.  \hfill$\square_{\scite{2.1}}$ 
\bigskip

\proclaim{\stag{2.5} Claim}  Assume that clause (f) of \scite{2.1}
fails, $d = d^* \in {\frak D}$ but clause (g) of \scite{2.1} holds
(equivalently $\langle {\text{\rm val\/}}_d(x):x \in X \rangle$ is not
constant).  \ub{Then} the following holds:
\mr
\item "{$\boxdot_1$}"  there are no $\bar s^0 \in {\text{\rm
conv\/}}(V^*_0(d^*)),\bar s^1 \in {\text{\rm conv\/}}(V^*_1(d^*))$
such that $(\frac{\bold n}{2} -1,\frac{\bold n}{2} -1)$ lie on
${\text{\rm conv\/}}\{\bar s^0,\bar s^1\}$ and for some $\ell \in
\{0,1\}$ this set, {\rm conv}$\{\bar s^0,\bar s^1\}$,
 contains an interior point of {\rm conv}$(V^*_\ell(d^*))$
\sn
\item "{$\boxdot_2$}"   the lines $L^*_0 = 
\{(\frac{\bold n}{2} -1,y):y \in \Bbb R\},L^*_1 = 
\{(x,\frac{\bold n}{2} -1):x \in \Bbb R\}$ divides the plane; and
${\text{\rm conv\/}}(V^*(d^*))$ is
{\roster
\itemitem{ $(i)$ }  included in one of the four closed half planes \ub{or}
\sn
\itemitem{ $(ii)$ }  is disjoint to at least one of the closed
quarters minus $\{(\frac{\bold n}{2} -1,\frac{\bold n}{2} -1)\}$. 
\endroster}
\ermn
\endproclaim
\bigskip

\remark{\stag{2.5A} Remark}  1) Recall $V^*(d^*) = 
V^*_0(d^*) \cup V^*_1(d^*)$ and $V^*_\ell(d^*) = 
\{\bar s - (\ell,1-\ell):\bar s \in V_\ell(d^*)\}$. \nl
2) So
\mr
\widestnumber\item{$(iii)$}
\item "{$(i)$}"  $(k_1,k_2) \in V^*_0(d^*) \Leftrightarrow (k_1,k_2) +
(0,1) \in V_0(d^*) \Leftrightarrow (k_1,k_2+1) \in
V_0(d^*)$
\nl
[Why?  By Definition \scite{1.10}(4).]
\sn
\item "{$(ii)$}"  $(k_2,k_1) \in V^*_1(d^*) 
\Leftrightarrow (k_2,k_1) + (1,0) \in V_1(d^*) 
\Leftrightarrow (k_2+1,k_1) \in V_1(d^*)$ \nl
(see \scite{1.11}(2)) \nl
hence
\sn
\item "{$(iii)$}"   $(k_1,k_2) \in V^*_0(d^*) \Leftrightarrow
(k_1,k_2 + 1) \in V_0(d^*) \Leftrightarrow (k_2 +1,k_1) \in V_1(d^*)
\Leftrightarrow (k_2,k_1) \in V^*_1(d^*)$. \nl
[Why? By the above $(i) + (ii)$ and \scite{1.11}(2).]
\endroster
\endremark
\bigskip

\demo{Proof}  Toward contradiction assume that $\boxdot_2$ or
$\boxdot_1$ in the claim fails.  So necessarily
\mr
\item "{$(*)_0$}"  $(\frac{\bold n}{2}-1,\frac{\bold n}{2} -1) \notin
V^*(d^*)$ \nl
[Why?  If it belongs to  $V^*_\ell(d^*)$ let $r^*_\ell =1,r^*_{1 - \ell}
= 0$ and we get clause (f) of \scite{2.1} which we are assuming fails]
\sn
\item "{$(*)'_0$}"  $(\frac{\bold n}{2}-1,\frac{\bold n}{2} -1) \notin
\text{ conv}(V^*_\ell(d^*))$ \nl
[Why?  As in the proof of $(*)_0$.]
\sn
\item "{$(*)_1$}"  $(\frac{\bold n}{2}-1,\frac{\bold n}{2}-1)$ belongs
to the convex hull of $V^*_0(d^*) \cup V^*_1(d^*)$ hence of
conv$(V^*_0(d^*)) \cup \text{ conv}(V^*_1(d^*))$ \nl
[Why?  Otherwise $\boxdot_1$ trivially holds; also there is a line $L$ through
$(\frac{\bold n}{2}-1,\frac{\bold n}{2}-1)$ such that $V^*(d^*)
\backslash L$ lie in one half plane of $L$, so easily clause (ii) of
$\boxdot_2$ holds hence $\boxdot_2$ holds recalling $(*)_0$.  
But we are assuming toward
contradiction that $\boxdot_1$ fails or $\boxdot_2$ fails.]
\ermn
Let $E = \{(\bar s_0,\bar s_1):\bar s_\ell \in \text{
conv}(V^*_\ell(d^*))$ for $\ell=0,1$ and 
$(\frac{\bold n}{2}-1,\frac{\bold n}{2}-1)$
belongs to the convex hull of $\{\bar s_0,\bar s_1\}\}$
\mr
\item "{$(*)_2$}"  $E \ne \emptyset$ \nl
[Why?  By $(*)_1$]
\sn
\item "{$(*)_3$}"  if $r_0,r_1 \in [0,1]_{\Bbb R},1 = r_0 + r_1,
(\frac{\bold n}{2}-1,\frac{\bold n}{2} -1) = r_0 \times \bar s_0 + r_1
\times \bar s_1$ and $\bar s_\ell \in \text{ conv}(V^*_\ell(d^*))$ for
$\ell =0,1$  \ub{then} $r_0 = r_1
= \frac 12$ \nl
[Why?  Otherwise clause (f) holds contradicting an assumption of \scite{2.5}.]
\sn
\item "{$(*)_4$}"  if $(\bar s_0,\bar s_1) \in E$ 
then $(\frac{\bold n}{2}-1,\frac{\bold n}{2}-1) = \frac 12(\bar s_0 +
\bar s_1)$ \nl
[Why?  By $(*)_3$ and the definition of $E$]
\sn
\item "{$(*)_5$}"  if $(\bar s_0,\bar s_1) \in E,\ell \in \{0,1\}$,
\ub{then} $\bar s_\ell$ is the unique member of conv$(V^*_\ell(d^*))$
which lies on the line through $\{\bar s_0,\bar s_1\}$.
\nl
[Why?  Otherwise let $\bar s'_\ell$ be a counterexample.  If
$(\frac{\bold n}{2}-1,\frac{\bold n}{2}-1) \in \text{ conv}\{\bar
s_\ell,\bar s'_\ell\}$ then it belongs to conv$(V^*_\ell(d^*))$ 
contradicting $(*)'_0$.  So letting
$\bar s'_{1-\ell} = \bar s_{1-\ell}$ we know that 
$(\frac{\bold n}{2}-1,\frac{\bold n}{2} -1) \in \text{ conv}\{\bar
s'_{1 - \ell},\bar s'_\ell\}$ hence by the definition of $E$ we get
$(\bar s'_0,\bar s'_1) \in E$ so by $(*)_4$ we deduce 
$\frac 12(\bar s'_0 + \bar s'_1) =
(\frac{\bold n}{2}-1,\frac{\bold n}{2}-1) = \frac 12(\bar s_0 + \bar
s_1)$ hence subtracting the two equations, $\bar s_{1 - \ell}$ is
cancelled and we get $\bar s'_\ell = \bar s_\ell$, contradiction]
\sn
\item "{$(*)_6$}"  $\boxdot_1$ holds (so by the assumption towards
contradiction $\boxdot_2$ fails). \nl
[Why?  Assume $\bar s_0,\bar s_1$ are as there hence (by the
definition of $E$), $(\bar s_0,\bar s_1) \in E$, now by $(*)_5$ the set (the
line through $\bar s_0,\bar s_1$) $\cap \text{ conv}(V^*_\ell(d^*))$ is
equal to $\{\bar s_\ell\}$.  So the line through $\bar s_0,\bar s_1$
cannot contain an interior point of 
conv$(V^*_\ell(d^*))$.]
\ermn
Easily (by $(*)_4$ and the definition of $E$):
\mr
\item "{$(*)_7$}"   $E_\ell := 
\{\bar s_\ell:(\bar s_0,\bar s_1) \in E\}$ is a convex subset
of conv$(V^*_\ell(d^*)) \subseteq \Bbb R^2$.
\ermn
Also
\mr
\item "{$(*)_8$}"  $(\frac{\bold n}{2} -1,\frac{\bold n}{2} -1) \notin
E_\ell$ \nl
[Why?  By $(*)'_0 + (*)_7$; also because if 
$(\frac{\bold n}{2} - 1,\frac{\bold n}{2} -1) \in
E_\ell$ then by $(*)_4$ we have $(\frac{\bold n}{2} -1,\frac{\bold
n}{2} -1) \in \text{ conv}(V^*_0(d^*)) \cap \text{ conv}(V^*_1(d^*))$
contradiction to $(*)'_0$.]
\endroster
\enddemo
\bn
Now we split the rest of the 
proof to three cases which by $(*)_2$ trivially exhausts all the
possibilities. \nl
\ub{Case 1}:  $E$ is not a singleton. 
\nl
This implies by $(*)_4$  that $E_\ell$ (defined in $(*)_7$) has at
 least two members, so by $(*)_7$ the set 
$V^*_\ell(d^*)$ is not a singleton for
$\ell=0,1$.  As $|E| \ge 2$ by $(*)_4$
clearly $|E_1| \ge 2$.  Also by $(*)_5$ if 
$\bar s_1 \in E_1$ so $\bar s_1 \ne (\frac{\bold n}{2} -1,\frac{\bold
n}{2} -1)$ by $(*)_8$, then $\bar s_1$ is
the unique member of conv$(V^*_1(d_1)) \cap$ (the line through $\bar
s_1, (\frac{\bold n}{2}-1,\frac{\bold n}{2}-1))$.  Also $E_1$ is
convex (by $(*)_7$) so necessarily
\mr
\item "{$(*)_9$}"   $E_1$ lies on a line $L_1$ to which by $(*)_5$, the point
$(\frac{\bold n}{2}-1,\frac{\bold n}{2}-1)$ does not belong.
\ermn
Let
\mr
\item "{$(*)_{10}$}"  $L_0$ is the line $\{(a_0,a_1):(a_1,a_0) \in
L_1\}$.
\ermn
As $E_1 \subseteq \text{ conv}(V^*_1(d^*)) \cap L_1$ is 
a convex set with $\ge 2$
members and $(*)_5$ it follows that conv$(V^*_1(d^*))$ is included 
in this line $L_1$ and as $V^*_0(d^*) = \{(k_2,k_1):(k_1,k_2) \in
V_1(d^*)\}$, (by clause $(iii)$ of Remark 
\scite{2.5A}(2) above) it follows that
conv$(V^*_0(d^*))$ is included in the line $L_0 = \{(a_0,a_1):(a_1,a_0) \in
L_1\}$ to which $(\frac{\bold n}{2} -1,\frac{\bold n}{2} -1)$ does not belong.

But $E_0 = \{\bar s_0:(\bar s_0,\bar s_1) \in E\}$ is necessarily an
interval of $L_0$ and by $(*)_4$ we have
\mr
\item "{$\odot_0$}"  $L_0 = \{(a_0,a_1):
2(\frac{\bold n}{2}-1,\frac{\bold n}{2} -1) - (a_0,a_1) \in L_1\}$.
\ermn
As $L_1$ is a line, for some reals $r_0,r_1,r_2$ we have
\mr
\item "{$\odot_1$}"  $L_1 = \{(a_0,a_1) \in 
\Bbb R^2:r_0 a_0 + r_1 a_1+r_2 = 0\}$
\ermn
and
\mr
\item "{$\odot_2$}"  $(r_0,r_1) \ne (0,0)$.
\ermn
Hence by the definition of $L_0$ in $(*)_{10}$ above we have
\mr
\item "{$\odot_3$}"  $L_0 = \{(a_0,a_1) \in \Bbb R^2:r_1a_0+r_0 a_1+r_2 = 0\}$
\ermn
and by $(*)_4$ the line $L_0$ includes the interval $\{2(\frac{\bold
n}{2} -1,\frac{\bold n}{2} -1) - \bar s_1:\bar s_1 \in E_1\}$ so

$$
L_0 = \{(a_0,a_1):(-r_0)a_0+(-r_1)a_1+r'_2 = 0\}
$$
\mn
where $r'_2 = 
2r_0(\frac{\bold n}{2}-1) + 2r_1(\frac{\bold n}{2} -1) + r_2$.

So for some $s \in \Bbb R$ we have $r_0 = -sr_1,r_1 = -sr_0,r_2 =
sr'_2$ but $(r_0,r_1) \ne (0,0)$
hence $s \in \{1,-1\}$ hence $r_0 \in \{r_1,-r_1\}$, so \wilog
\, $r_0=1,r_1 \in \{1,-1\}$.
\bn
\ub{Subcase 1A}:  $r_1 = -1$.

So $d^*\{x,y\} = y \Rightarrow (\text{val}_{d^*}(x),\text{val}_{d^*}(y)) \in
V_1(d^*) \Rightarrow (\text{val}_{d^*}(x)-1,
\text{ val}_{d^*}(y)) \in V^*_1(d^*) \Rightarrow
(\text{val}_{d^*}(x)-1,\text{val}_{d^*}(y)) \in L_1 
\Rightarrow \text{ val}_{d^*}(x) 
- \text{ val}_{d^*}(y) = -r_2 +1$, i.e., 
is constant, is the same for any such pair $(x,y)$.  
If the directed graph Tor$(d^*) = (X,\{(u,v):d^*\{u,v\} = v\})$ 
contains no cycle, or just no cycle of length 3, then for some
list $\{x_\ell:\ell < \bold n\}$ of $X$, we have
$d^*\{x_{\ell_1},x_{\ell_2}\} = x_{\text{max}\{\ell_1,\ell_2\}}$ for
$\ell_1 \ne \ell_2 < \bold n$.  This implies $V_0(d^*) =
\{(\ell_1,\ell_2):\ell_1 < \ell_2 < \bold n\}$, easy
contradiction to $(*)'_0$.  So the directed graph Tor$(d^*) =
(X,\{(u,v):d^*\{u,v\} = v\})$ necessarily contains 
a cycle, so necessarily $-r_2+1=0$.
Recall that when ${\frak D} \subseteq {\frak C}^{\text{full}}$, the
graph is connected so the val$_{d^*}(x)$ is the same 
for all $x \in X$ hence is
necessarily $(\frac{\bold n}{2}-1)$, which is not an ``allowable" case, in
particular, contradict clause (g) of \scite{2.1} which we are assuming.
\bn
\ub{Subcase 1B}:  $r_1 = 1$.

Clearly $x \ne y \in X \and d^*\{x,y\} = y \Rightarrow
(\text{val}_{d^*}(x),\text{val}_{d^*}(y)) \in V_1(d^*) \Rightarrow
(\text{val}_{d^*}(x)-1,\text{val}_{d^*}(y)) \in V^*_1(d^*) \Rightarrow
(\text{val}_{d^*}(x)-1,\text{val}_{d^*}(y)) \in L_1 \Rightarrow 
\text{ val}_{d^*}(x) + \text{ val}_{d^*}(y) = -r_2+1$.  As $\bold n
\ge 3$ and recall that ${\frak D} \subseteq {\frak C}^{\text{full}}$, so
there are distinct $x_0,x_1,x_2 \in X$ so val$_{d^*}
(x_{\ell_1}) + \text{ val}_{d^*}(x_{\ell_2}) = -r_2+1$ for 
$\{\ell_1,\ell_2\} \in \{\{0,1\},\{0,2\},\{1,2\}\}$, the order is not
important as $r_1 = r_0$ hence val$_{d^*}(x_1)$, val$_{d^*}(x_2)$
 are equal and $-r_2+1$ is twice their
value.  So for $y \in X \backslash \{x_1\}$, we have val$_{d^*}(x_1) +
\text{ val}_{d^*}(y) = -r_2+1$ so val$_{d^*}(y) = 
\text{ val}_{d^*}(x_2)$, so we are done as in case 1A.
\bn
\ub{Case 2}:  $E = \{(\bar s^*_0,\bar s^*_1)\}$ and $\bar s^*_0 \ne
\bar s^*_1$.

Let $L$ be the line through $\{\bar s^*_0,\bar s^*_1\}$ and let the
real $r_0,r_1,r_2$ be such that $L= 
\{(a_0,a_1):r_0a_0+r_1a_1+r_2=0$ and $(r_0,r_1) \ne (0,0)\}$.

So by $(*)_5$ the set 
conv$(V^*_\ell(d^*))$ intersect $L$ in the singleton $\{\bar
s^*_\ell\}$
\mr
\item "{$(*)_{10}$}"  no one (closed) half plane for the 
line $L$ contains $V^*_0(d^*) \cup V^*_1(d^*)$.
\ermn
[Why?  As then $\boxdot_2$ of \scite{2.5} holds (if $L$ is not
parallel to the $x$-axis and the $y$-axis (i.e., $r_0,r_1 \ne 0$) then
$\boxdot_2(ii)$ holds, otherwise $\boxdot_2(i)$ holds); so by $(*)_6$
this is against our assumption toward contradiction.]
\mr
\item "{$(*)_{11}$}"  There are $\bar k_0 \in V^*_0(d^*)
\backslash \{\bar s^*_0\}$ and $\bar k_1 \in V^*_1(d^*) \backslash
\{\bar s^*_1\}$ such that they are outside $L$ in different sides.
\ermn
[Why?  First, if there are $\ell \in \{0,1\}$ and $\bar k',\bar k'' \in
V^*_\ell(d^*) \backslash L$ on different sides of $L$ then also
$V^*_{1-\ell}(d^*)$ has a member outside $L$ (by clause $(iii)$ of
\scite{2.5A}(2) above and $(*)_5$) and call it $\bar k_{1 - \ell}$, so the
choice $\bar k_\ell = \bar k'$ or the choice $\bar k_\ell = \bar k''$
is as required.  Second, if there is no such $\ell \in \{0,1\}$ by 
$(*)_{10}$ we can choose $\bar k_\ell \in V^*_\ell(d^*)
\backslash L$ for $\ell=0,1$ and by the same reason 
$\bar k_0,\bar k_1$ are as
required.]

As $(\frac{\bold n}{2}-1,\frac{\bold n}{2}-1)$ lie in the open interval
spanned by $\bar s^*_0$ and $\bar s^*_1$, necessarily
$(\frac{\bold n}{2}-1,\frac{\bold n}{2}-1)$ is an interior point of 
conv$\{\bar s^*_0,\bar k_0,\bar s^*_1,\bar k_1\}$, easy contradiction to the
case assumption.
\bn
\ub{Case 3}:  $E = \{(\bar s^*_0,\bar s^*_1)\}$ and $\bar s^*_0 = \bar s^*_1$.

So by $(*)_4$ clearly 
$\bar s^*_\ell = (\frac{\bold n}{2} -1,\frac{\bold n}{2}-1)$, but
this contradicts $(*)'_0$ above. \hfill$\square_{\scite{2.5}}$
\bigskip

\proclaim{\stag{2.6} Claim}  In \scite{2.5}, clause (i) of $\boxdot_2$ 
is impossible.
\endproclaim
\bigskip

\demo{Proof}  So toward contradiction assume clause (i) of $\boxdot_2$ holds.
 We know that $\langle \text{val}_{d^*}(x):
x \in X \rangle$ is not constant (as we assume clause (g) of
\scite{2.1} holds).  As the average of val$_{d^*}(x),x \in X$ is
$\frac{\bold n-1}{2} = \frac{\bold n}{2}- \frac 12$ clearly
\mr
\item "{$(*)_1$}"  for some points $x \in X$ we have val$_{d^*}(x)$ is 
$< \frac{\bold n-1}{2} = \frac{\bold n}{2} - \frac 12$ and for some
point $x \in X$ we have val$_{d^*}(x)$ is $> \frac{\bold n}{2} - \frac 12$.
\ermn
The assumption (i.e. (i) of $\boxdot_2$ of \scite{2.5})
leaves us with four possibilities, so we have 4 cases (according to
which half plane).  
\enddemo
\bn
\ub{Case 1}:  For no $(k_0,k_1) \in V^*(d^*)$ do we have $k_0 >
\frac{\bold n}{2}-1$.

It follows that by clause (i) of \scite{2.5A}(2)

$$
(k_0,k_1) \in V_0(d^*) \Rightarrow (k_0,k_1-1) \in V^*_0(d^*)
\subseteq V^*(d^*) \Rightarrow k_0 \le \frac{\bold n}{2}-1.
$$
\mn
So if $x \in X$ and for some $y \in X \backslash \{x\}$ we have
$d^*\{x,y\} = x$ (this means just that, val$_{d^*}(x) < \bold n-1)$ then
val$_{d^*}(x) \le \frac{\bold n}{2}-1$, so
\mr
\item "{$(*)_2$}"  if $x \in X$ and val$_{d^*}(x) < \bold n-1$ \ub{then}
val$_{d^*}(x) \le \frac{\bold n}{2}-1$.
\ermn
We shall show that this is impossible (this helps also in case 3).
There can be at most one $x \in X$ with val$_{d^*}(x) = \bold n -1$; if
there is none then we have:

if $x \in X$ then val$_{d^*}(x) \le \frac{\bold n}{2} - 1$.

But the average valency is $\frac{\bold n -1}{2}$ which is $>
\frac{\bold n}{2} -1$, contradiction.  So 
there is $x^* \in X$ such that val$_{d^*}(x^*) = \bold n-1$,
of course, it is unique.  Now Tor$^-(d^*) := 
(X \backslash \{x^*\},\{(y,z):y \ne z
\in X \backslash \{x^*\},d^*\{y,z\} = z\})$ is a directed graph 
with $\bold n-1$ points  and every $y \in X \backslash \{x^*\}$ has
the same out-valency in Tor$(d^*)$ and in Tor$^-(d^*)$, hence
each $y \in X \backslash \{x^*\}$ has (in Tor$(d^*)$ and in
Tor$^-(d^*)$) out-valency $\le \frac{\bold n}{2}-1 = 
\frac{(\bold n-1)-1}{2}$, so necessarily $\bold n$ is even and every
node in Tor$^-(d^*)$ has out-valency exactly $\frac{(\bold n
-1)-1}{2} = \frac{\bold n}{2} - 1$; as 
$\bold n \ge 3$ we can choose $y \ne z \in
X \backslash \{x^*\}$ and \wilog \, $d^*(y,z) = z$.  Now
$(\frac{\bold n}{2}-1,\frac{\bold
n}{2}-1),(\bold n-1,\frac{\bold n}{2}-1) \in V_1(d^*)$ as witnessed by
the pairs $(y,z),(x^*,y)$ respectively, hence
$(\frac{\bold n}{2}-2,\frac{\bold n}{2}-1),(\bold n-2,\frac{\bold
n}{2}-1) \in V^*_1(d)$, again by Definition \scite{1.10}(4). 
Hence (as $\bold n \ge 3$ so $\bold n -2 \ge \frac{\bold n}{2} -1$) we
have $(\frac{\bold n}{2} -1,\frac{\bold n}{2}- 1) \in 
\text{ conv}(V^*_1(d^*))$, so $r^*_0 = 0$ given 
contradiction to ``clause (f) of
\scite{2.1} fails" assumed in \scite{2.5}.
\bn
\ub{Case 2}:  For no $(k_0,k_1) \in V^*(d^*)$ do we have $k_0 <
\frac{\bold n}{2}-1$.  \nl
By clause (i) of \scite{2.5A}(2) it follows that

$$
(k_0,k_1) \in V_0(d^*) \Rightarrow (k_0,k_1-1) \in V^*_0(d^*)
\subseteq V^*(d^*) \Rightarrow k_0 \ge \frac{\bold n}{2}-1.
$$
\mn
So if $x \in X$ and for some $y \in X \backslash \{x\}$ we have
$d^*\{x,y\} = x$ (equivalently val$_{d^*}(x) < \bold n -1)$
then $(\text{val}_{d^*}(x),\text{val}_{d^*}(y)) \in V_0(d^*)
\Rightarrow (\text{val}_{d^*}(x),\text{val}_{d^*}(y)-1) \in V^*_0(d^*)$ hence
val$_{d^*}(x) \ge \frac{\bold n}{2}-1$, but $\bold n -1 \ge
\frac{\bold n}{2} -1$ so in any case
\mr
\item "{$(*)_3$}"  if $x \in X$ then val$_{d^*}(x) \ge \frac{\bold
n}{2}-1$.
\ermn
For the rest of the proof of case 2, we shall use $(*)_3$.  This
serves us also in case 4.
So $x \in X \Rightarrow \text{ val}_{d^*}(x) \ge \frac{\bold n}{2}-1$.

If $\bold n$ is odd we have $x \in X \Rightarrow \text{ val}_{d^*}(x)
\ge \frac{\bold n}{2} - \frac 12 = \frac{\bold n -1}{2}$, impossible
by $(*)_1$ so $\bold n$ is even.  Let $k = \frac{\bold n}{2}-1$.  The
average val$_{d^*}(x)$ is necessarily $k + \frac 12$ hence $Y =: \{x \in
X:\text{val}_{d^*}(x) \le k$ (equivalently $=k)\}$ has at least 
$k+1 = \frac{\bold n}{2}$ members.  If $x \in
X$,val$_{d^*}(x) = k+1$ then $x \notin Y$ and $|\{y:d\{x,y\} = y\}| = k+1 =
\frac{\bold n}{2} > |X \backslash (Y \cup \{x\})|$ 
so there is $y \in Y$ such that $d\{x,y\} = y$, hence $(k,k) =
(\text{val}_{d^*}(x)-1$, val$_{d^*}(y)) \in V^*_1(d^*)$ so clause (f) of
\scite{2.1} holds with $r^*_0 =0$, contradiction. So
\mr
\item "{$(*)_4$}"  $x \in X \Rightarrow \text{ val}_{d^*}(x) \ne k+1$.
\ermn
Now $|Y| = \bold n$ is impossible by $(*)_1$.  Also if $|Y| = \bold
n-1$ let $x^*$ be the unique element of $X$ outside $Y$ so in the
tournament Tor$^- := (Y,\{(x,y):x \ne y$ are from 
$Y$ and $d^*\{x,y\} = y\})$
each $x$ has out-valency $\le$ val$_{d^*}(x) = k = \frac{(|Y|-1)}{2}$,
but this is the average so equality holds.  Now if 
$x \in Y$ then $x$ has out-valency
$k$ in Tor$^-(d^*)$ and has out-valency $k$ in Tor$(d^*)$ hence
$d^*\{x,x^*\} \ne x^*$ hence val$_{d^*}(x^*) = \bold n-1$
 and we get contradiction as in Case 1.  Hence
\mr
\item "{$(*)_5$}"  $|Y| \le \bold n -2$.
\ermn
Clearly we can find $x_1 \in Y$ such that $|\{y \in Y:y \ne x_1,d^*\{x_1,y\} =
y\}| \le \frac{|Y|-1}{2}$ (as if we average this number on the $x_1
\in Y$ we get $\frac{|Y|-1}{2}$) but 
$\frac{|Y|-1}{2} \le \frac{\bold n}{2} - \frac 32 = 
k - \frac 12 < |\{y \in X:d\{x_1,y\} = y\}|$ 
hence there is $x_2 \in X \backslash Y$ such that
$d^*\{x_1,x_2\} = x_2$.  Now
let $m = \text{ val}_{d^*}(x_2)$ so $m > k$ as $x_2 \notin Y$ and $m
\ne k+1$ by $(*)_4$ hence $m > k+1$ and $(x_1,x_2)$ witness $(k-1,m) \in
V^*_1(d^*)$.  As $\bold n \ge 3$ and (see the paragraph before $(*)_4$) 
$|Y| \ge \frac{\bold n}{2}$
obviously $|Y| \ge 2$ hence (as any pair of $y_1 \ne y_2$ from $Y$
witness) also $(k-1,k) \in V^*_1(d^*)$.  As
$|Y| \ge \frac{\bold n}{2} = k+1$, val$_{d^*}(x_2) > k+1 \ge \bold n-|Y|$
easily there is $x_3 \in Y$ such that $d\{x_2,x_3\} = x_3$ hence
$(x_2,x_3)$ witness $(m-1,k) \in V^*_1(d^*)$.  Now 
$(\frac{\bold n}{2} -1,\frac{\bold n}{2} -1) = (k,k)
\in \text{ conv}\{(k-1,k),(m-1,k)\})$ recalling $m > k+1$.
But $(k-1,k),(m-1,k)$ belong to $V^*_1(d^*)$ hence $(\frac{\bold n}{2}
-1,\frac{\bold n}{2} -1) \in \text{ conv}(V^*_1(d^*))$,
contradiction to ``not clause (f) of \scite{2.1}" with $r^*_0=0$.  
\bn
\ub{Case 3}: For no $(k_0,k_1) \in V^*(d^*)$ do we have $k_1 >
\frac{\bold n}{2} -1$.

So by clause (ii) of \scite{2.5A} it follows that

$$
(k_0,k_1) \in V_1(d^*) \Rightarrow (k_0 -1,k_1) \in V^*_1(d^*)
\subseteq V^*(d^*) \Rightarrow k_1 \le \frac{\bold n}{2} -1.
$$
\mn
So if $y \in X$ and for some $x \in X \backslash \{y\}$ we have
$d\{x,y\} = y$ then $(\text{val}_{d^*}(x)-1,\text{val}_{d^*}(y)) \in
V^*_1(d^*)$ hence val$_{d^*}(y) \le \frac{\bold n}{2} -1$.  But there is
such $x$ iff val$_{d^*}(y) \ne \bold n -1$, that is
\mr
\item "{$(*)_6$}"  if $y \in X$ and val$_{d^*}(y) \ne \bold n -1$ then
val$_{d^*}(y) \le \frac{\bold n}{2} -1$.
\ermn
We continue as in Case 1, (after $(*)_2$ which uses only $(*)_2$ or 
dualize see \scite{1.11}(1)).
\bn
\ub{Case 4}:  For no $(k_0,k_1) \in V^*(d^*)$ do we have $k_1 <
\frac{\bold n}{2} -1$.  

So $(k_0,k_1) \in V_1(d^*) \Rightarrow (k_0-1,k_1) \in V^*_1(d^*)
\subseteq V^*(d^*) \Rightarrow k_1 \ge \frac{\bold n}{2} -1$.  So if
$y \in X$ and for some $x \in X \backslash \{y\}$ we have $d\{x,y\} =
y$ then $(\text{val}_{d^*}(x)-1,\text{val}_{d^*}(y)) \in V^*_1(d^*)
\Rightarrow \text{ val}_{d^*}(y) \ge \frac{\bold n}{2} -1$.  So if 
val$_{d^*}(y) < \bold n-1$ then there is such $x$ hence
val$_{d^*}(y) \ge \frac{\bold n}{2} -1$, but 
if val$_{d^*}(y) \ge \bold n -1$ we get the same
conclusion, so
\mr
\item "{$(*)_7$}"  val$_{d^*}(y) \ge \frac{\bold n}{2} -1$
\ermn
and we can continue as in case 2 after $(*)_3$.  \hfill$\square_{\scite{2.6}}$
\bigskip

\proclaim{\stag{2.7} Claim}  In \scite{2.5}, clause (ii) of $\boxdot_2$
is impossible.
\endproclaim
\bigskip

\demo{Proof}  Note that as we are assuming the failure of
clause (f) of \scite{2.1}
\mr
\item "{$(*)_0$}"  $(\frac{\bold n}{2} -1,\frac{\bold n}{2} -1) \notin
V^*(d^*)$.
\ermn
Again we have four cases.
\sn
\ub{Case 1}:  If $a_0 \ge \frac{\bold n}{2} -1,a_1 \ge \frac{\bold
n}{2} -1$ but $(a_0,a_1) \ne (\frac{\bold n}{2} -1,\frac{\bold n}{2} -1)$
then $(a_0,a_1) \notin \text{ conv}(V^*(d^*))$. \nl
So
\mr
\item "{$(*)_1$}"  for at most one $x \in X$ we have val$_{d^*}(x) \ge
\frac{\bold n}{2}$.\nl
[Why?  If $x \ne y \in X$ and val$_{d^*}(x) \ge \frac{\bold n}{2}$,
val$_{d^*}(y) \ge \frac{\bold n}{2}$ \ub{then}
$(\text{val}_{d^*}(x)-1,\text{val}_{d^*}(y)) \in V^*_1(d^*) \subseteq
V^*(d^*)$ or $(\text{val}_{d^*}(x),\text{val}_{d^*}(y)-1) \in
V^*_0(d^*) \subseteq V^*(d^*)$, a contradiction to the case assumption
in both cases.]
\ermn
If there is no $x \in X$ with val$_{d^*}(x) \ge \frac{\bold n}{2}$ 
then $x \in X \Rightarrow \text{ val}_{d^*}(x) <
\frac{\bold n}{2}$ and so $x \in X \Rightarrow \text{ val}_{d^*}(x)
\le \frac{\bold n-1}{2}$ but this is the average valency, so always equality
holds, contradicting an assumption of \scite{2.5}.
\nl
So assume
\mr
\item "{$(*)_2$}"  $x_0 \in X$, val$_{d^*}(x_0) \ge \frac{\bold
n}{2}$.
\ermn
Now
\mr
\item "{$(*)_3$}"  if $y \in X \backslash \{x_0\}$ and $d^*\{x_0,y\} =
y$ \ub{then} $\text{val}_{d^*}(y) < \frac{\bold n}{2} -1$ 
(hence $\le \frac{\bold n}{2} - \frac 32$). 
\nl
[Why?  As $d^*\{x_0,y\} = y$ then $(\text{val}_{d^*}(x_0)-1,
\text{val}_{d^*} (y)) \in V^*_1(d^*) \subseteq V^*(d^*)$, now
$\text{val}_{d^*}(x_0) -1 \ge \frac{\bold n}{2} -1$ hence by the case
assumption $+ (*)_0$ we have $\text{val}_{d^*}(y) < \frac{\bold n}{2}
-1$.]
\sn
\item "{$(*)_4$}"  if $y \in X \backslash \{x_0\}$ and $d^*\{x_0,y\} =
x_0$ \ub{then} $\frac{\bold n}{2} -1 > |\{z \in X \backslash
\{x_0,y\}:d^*\{y,z\} = z\}|$.
\nl
[Why?  As $d^*\{x_0,y\} = x_0$ clearly $(\text{val}_{d^*}(x_0),
\text{val}_{d^*}(y)-1) \in V^*_0(d^*) \subseteq V^*(d^*),
\text{val}_{d^*}(x_0) \ge \frac{\bold n}{2} > \frac{\bold n}{2}-1$, 
hence by the case assumption $\text{val}_{d^*}(y) -1 < \frac{\bold
n}{2} -1$ so
val$_{d^*}(y) < \frac{\bold n}{2}$, i.e. $\frac{\bold n}{2} > |\{z \in
X \backslash \{y\}:d^*\{y,z\} = z\}|$ and as $d^*\{x_0,y\} = x_0$ this
gives the desired inequality.]
\ermn
So letting $Y = X \backslash \{x_0\}$ we have
$(Y,\{\{y,z\},y \ne z \in Y,d^*(y,z) = z\})$ is a tournament
satisfying each node has out-valency $\le \frac{\bold
n-3}{2} < \frac{(\bold n-1)-1}{2} = \frac{|Y|-1}{2}$ (why? by
$(*)_3+(*)_4$), contradiction.
\enddemo
\bn
\ub{Case 2}:  If $a_1 \le \frac{\bold n}{2}-1$ and $a_2 \le
\frac{\bold n}{2}-1$ and $(a_1,a_2) \ne (\frac{\bold n}{2}
-1,\frac{\bold n}{2} -1)$ \ub{then} $(a_1,a_2) \notin \text{ conv}(V^*(d^*))$.
\nl
Clearly, as above in the proof of $(*)_1$
\mr
\item "{$(*)'_1$}"  there is at most one $x \in X$ with val$_{d^*}(x)
\le \frac{\bold n}{2}-1$.
\ermn
If there is none then $x \in X \Rightarrow \text{ val}_{d^*}(x) \ge
\frac{\bold n}{2}-1 + \frac 12 = \frac{\bold n-1}{2}$, so considering the
average of val$_{d^*}(y)$ equality always holds so clause (g) of
\scite{2.1} fails contradicting an assumption of \scite{2.5}.  So
assume
\mr
\item "{$(*)'_2$}"  $x_0 \in X$, val$_{d^*}(x_0) \le \frac{\bold n}{2}-1$
\ermn
and by $(*)'_1 + (*)'_2$ clearly
\mr
\item "{$(*)'_3$}"  if $y \in X \backslash \{x_0\}$
then val$_{d^*}(y) > (\frac{\bold n}{2}-1)  = \frac{\bold n -2}{2}$ so
val$_{d^*}(y) \ge \frac{\bold n-1}{2}$.
\ermn
The directed graph $\bold G = (X \backslash \{x_0\},\{(y,z):d\{y,z\} = z\})$
has $\bold n-1$ nodes and let $Y_0 = \{y \in X:y \ne x_0$ and
$d^*\{y,x_0\}=x_0\}$ and $Y_1 = \{y \in X:y \ne x_0$ and $d^*\{x_0,y\} =
y\}$.

Clearly $Y_0,Y_1,\{x_0\}$ is a partition of $X$, so $Y_0,Y_1$ is a
partition of the set of nodes in $\bold G$.  Also
\mr
\widestnumber\item{$(iii)$}
\item "{$\boxdot$}"  $(i) \quad y \in Y_0 \Rightarrow 
d^*\{x,y\} = x_0 \Rightarrow (\text{val}_{d^*}(x_0),
\text{val}_{d^*}(y)-1) \in V^*_0(d^*) \subseteq$
\nl

\hskip25pt $V^*(d^*) \Rightarrow$ (by the case assumption $+ (*)'_2 + (*)_0)
\text{ val}_{d^*}(y)-1 > \frac{\bold n}{2} -1 \Rightarrow \text{
val}_{d^*}(y) >$
\nl

\hskip25pt  $\frac{\bold n}{2} \Rightarrow$ 
the valency of $y$ in $\bold G$ is 
$> \frac{\bold n}{2} -1 \Rightarrow$ 
\nl

\hskip25pt the valency of $y$ in $\bold G$
is $\ge \frac{\bold n -1}{2}$
\sn
\item "{${{}}$}"  $(ii) \quad y \in Y_1 \Rightarrow d^*\{x_0,y\} = y$
(by $(*)_3$) $\Rightarrow \text{ val}_{d^*}(y) \ge 
\frac{\bold n -1}{2} \Rightarrow$ the valency 
\nl

\hskip25pt  of $y$ in $\bold G$ is
$\ge \frac{\bold n -1}{2} -0 = \frac {\bold n -1}{2}$.
\ermn
So every node in $\bold G$ has out-valency (in
$\bold G$) at least $\frac{\bold n -1}{2}$, a contradiction as the
average out-valency is $\frac{\bold n-2}{2}$.
\bn
\ub{Case 3}:  If $a_1 \ge \frac{\bold n}{2}-1$ and $a_2 \le
\frac{\bold n}{2}-1$ and $(a_1,a_2) \ne (\frac{\bold n}{2}
-1,\frac{\bold n}{2} -1)$, \ub{then} $(a_1,a_2) \notin \text{ conv}(V^*(d^*))$.
\nl
So (as in the proof of $(*)_1$ using $(*)_0$)
\mr
\item "{$\bigodot_1$}"  there cannot be $x_0,x_1 \in X$ such that
val$_{d^*}(x_0) \ge \frac{\bold n}{2}$ and val$_{d^*}(x_1) \le
\frac{\bold n}{2}-1$ (the $x_0 \ne x_1$ follows)
\ermn
so one of the following two sub-cases hold.
\bn
\ub{Subcase 3A}:  $x \in X \Rightarrow \text{ val}_{d^*}(x) <
\frac{\bold n}{2}$.

So $x \in X \Rightarrow \text{ val}_{d^*}(x) \le \frac{\bold n-1}{2}$
and (looking at average valency) equality holds, contradicting clause
(g) of \scite{2.1} which we are assuming.
\bn
\ub{Subcase 3B}:  $x \in X \Rightarrow \text{ val}_{d^*}(x) >
\frac{\bold n}{2}-1$.

So $x \in X \Rightarrow \text{ val}_{d^*}(x) \ge \frac{\bold n-1}{2}$,
and we finish as above.
\bn
\ub{Case 4}: If $a_1 \le \frac{\bold n}{2} -1$ and $a_2 
\ge \frac{\bold n}{2}-1$ and $(a_1,a_2) \ne (\frac{\bold n}{2}
-1,\frac{\bold n}{2} -1)$ then $(a_1,a_2) \notin \text{ conv}(V(d^*))$.
\nl
Similar to case 3 (or dualize the situation by \scite{1.11}(1)).
\hfill$\square_{\scite{2.7}}$
\newpage

\head {\S3 Balanced choice functions} \endhead  \resetall \sectno=3
 \spuriousreset
\bigskip

Here we analyze the case clause (g) of \scite{2.1} fail and give a
complete answer (and show the equivalence of relatives of ``balance".
\bigskip

\definition{\stag{b.1} Definition}  1) $c \in {\frak C}$ is called
\ub{balanced} if $x \in X \Rightarrow \text{ val}_c(x) = (\bold n-1)/2$,
let ${\frak C}^{\text{bl}} = \{c \in {\frak C}:c$ is balanced$\}$. \nl
2) $\bar t \in \text{ pr}({\frak C})$ is called balanced \ub{if} $x
\in X \Rightarrow \Sigma\{t_{x,y}:y \in X \backslash \{x\}\} =
(\bold n-1)/2$.  Let $\text{pr}^{\text{bl}}({\frak C})$ be the set of balanced
$\bar t \in \text{ pr}({\frak C})$. 
\nl
2A) $\bar t = \text{ pr}({\frak C})$ is super-balanced if $t_{x,y} =
\frac 12$ for $x \ne y \in X$.
\nl
3) We say $c \in {\frak C}$ is pseudo-balance iff every edge of
   Tor$(c)$ belongs to a directed cycle, see Definition
   \scite{0.2}(2).
\nl
3A) $c \in {\frak C}$ is called partition$^+$-balanced when: if
$\emptyset \ne Y \subsetneqq X$ then for some $x \in X \backslash Y$
and $y \in Y$ we have $c\{x,y\} = y;
c \in {\frak C}$ is called partition-balanced when
``if $\emptyset \subsetneqq Y \subsetneqq X$" \ub{then}: for some $x \in
X \backslash Y,y \in Y$ we have $c\{x,y\} = y$ \ub{iff} for some $x
\in X \backslash Y,y \in Y$ we have $c(x,y) = x$.
\nl
3B) $c \in {\frak C}$ is called weight-balanced \ub{when}: for some
balanced $\bar t \in \text{ pr}({\frak C})$ we have

$$
c\{x,y\} = y \Leftrightarrow t_{x,y} > \frac{1}{2}.
$$
\mn
4) We call ${\Cal D} \subseteq {\frak C}$ balanced \ub{if} every $c
\in {\Cal D}$ is balanced, similarly $T \subseteq \text{ pr}({\frak
C})$ is called balanced if every $\bar t \in T$ is.  Similarly for the
other properties.
\nl
5) If $x,y,z \in X$ are distinct, let $\bar t = \bar t^{<x,y,z>}$
be defined by: 
\smallskip

$t_{u,v}$ is $1$ if $(u,v) \in \{(x,y),(y,z),(z,x)\}$ \nl
\smallskip

$t_{u,v}$ is $0$ if $(u,v) \in \{(y,x),(z,y),(x,z)\}$ \nl
\smallskip

$t_{u,v}$ is $\frac 12$ if otherwise. 
\sn
6) For $k \ge 3$ a sequence 
$\bar x = (x_0,\dotsc,x_{k-1})$ with $x_\ell \in X$ and no repetitions
and $a \in [0,1]_{\Bbb R}$ let $\bar t = 
\bar t_{\bar x,a} \in \text{ pr}({\frak C})$ be 
defined by $t_{x_i,x_j} = a,t_{x_j,x_i} = 1-a$ if
$j = i+1$ mod $k$ and $t_{x,y} = \frac 12$ for $x \ne y \in X$
otherwise.  If $a=1$ we may omit it.
\nl
7) Let $c_* \in {\frak C}$ be the empty
function.  We call ${\Cal D}$ trivial if ${\Cal D} = \{c_*\}$ 
or ${\Cal D} = \emptyset$.
\enddefinition
\bn
\margintag{b.2}\ub{\stag{b.2} Fact}  0) pr$^{\text{bl}}({\frak C})$ is a convex subset
of pr$({\frak C})$ and it is preserved by the permutations of
pr$({\frak C})$ induced by permutations of $X$.
\nl
1) If $c \in {\frak C}^{\text{bl}}$ then $\bar t[c]$
belongs to $\text{pr}^{\text{bl}}({\frak C})$ but is not necessarily
super-balanced; if $\bar t \in \text{ pr}({\frak C})$ is balanced 
and $c = \text{ maj}(\bar t)$,
\ub{then} $c$ is pseudo-balanced.
\nl
2) If $c \in {\frak C}^{\text{full}}$ is weight-balanced \ub{then} it
 is partition$^+$-balanced.
\nl
2A) If $c \in {\frak C}$ is weight-balanced \ub{then} it is partition-balanced
(similar to 
$(*)$ of \scite{1.9}; note that not every $c \in {\frak C}$ is
pseudo-balanced and even some 
$c \in {\frak C}^{\text{full}}$ is not pseudo-balanced). 
\nl
3) If $c \in {\frak C}^{\text{full}}$ \ub{then} $c$ is
   partition$^+$-balaned iff $c$ is partition-balanced...
\nl
4) If $c \in {\frak C}$ then maj$(\bar t[c]) = c$. 
If ${\Cal D} \subseteq {\frak C}$ is balanced, \ub{then} pr-cl$({\Cal D})$ is
balanced hence every member of maj-cl$({\Cal D})$ is pseudo-balanced.
\bigskip

\demo{Proof}  0), 1), 3), 4) Check.
\nl
2) By (2A) and 3).
\nl
2A) Toward contradiction assume $(Y,x,y)$ is a counterexample so
$\emptyset \subsetneqq Y \subsetneqq X,x \in X \backslash Y,y \in Y$
and $c\{x,y\} = y$ and there are no $u \in Y,v \in X \backslash Y$
such that $c\{u,v\} = v$.  Let $Y_0 = X \backslash Y,Y_1=X$.

As $c$ is pseudo-balanced, there is $\bar t \in \text{
pr}^{\text{bl}}({\frak C})$ such that $c = \text{ maj}(c)$.  So
necessarily $t_{x,y} > \frac 12$ but $u \in Y_0,v \in Y_1 \Rightarrow
t_{v,u} \le \frac 12 \Rightarrow t_{u,v} \ge \frac 12$.

Let $Y_0 = X \backslash Y,Y_1 = Y$, so $\Sigma\{t_{u,v} - \frac 12:u
\in Y_0,v \in Y\} \ge t_{x,y} - \frac 12 > 0$.

For $u \in Y_0$ let $s^0_u = \Sigma\{t_{u,v} - \frac 12:v \in Y_0
\backslash \{u\}\},s^1_u = \Sigma\{t_{u,v} - \frac 12:v \in Y_1\}$ so
$s^0_u + s^1_u = 0$, hence $0 = \Sigma\{s^0_u + s^1_u:u \in Y_0\} =
\Sigma\{s^0_u:u \in Y_0\} + \Sigma\{s^1_u:u \in Y_0\}$ but by the
previous sentence the second summand is positive.  Hence
$\Sigma\{s^0_u:u \in y_0\}$ is negative, but it is zero because $u_1 \ne
u_2 \in y \Rightarrow (t_{u_1,u_2} - \frac 12) + (t_{u_2,u_1} - \frac
12) =0$. 
\enddemo
\bigskip

\proclaim{\stag{3.1A} Claim}  If ${\frak D} \subseteq 
{\frak C}^{\text{full}}$ is non-empty, symmetric and 
not balanced, \ub{then} ${\text{\rm maj-cl\/}}({\frak D}) = {\frak C}$.
\endproclaim
\bigskip

\demo{Proof}  Choose $d^* \in {\frak D}$ which is not balanced, and
let ${\frak D}' =: \text{ sym-cl}(\{d^*\})$, so ${\frak D}'$ is as in
\scite{2.1} and it satisfies clause (g) there hence it satisfies clause
$(a)'$ there.  This means that maj-cl$({\frak D}') = {\frak C}$ but
${\frak D}' \subseteq {\frak D} \subseteq {\frak C}$ hence ${\frak C}
= \text{ maj-cl}({\frak D}') \subseteq \text{ maj-cl}({\frak D})
\subseteq {\frak C}$ so we are done.  \hfill$\square_{\scite{3.1A}}$
\enddemo
\bn
\margintag{b.2A}\ub{\stag{b.2A} Fact}:  1) If $c \in {\frak C}^{\text{bl}}$ or just
$c \in {\frak C}$ is weight-balanced \ub{then} $\bold c$ is
partition-balanced and is pseudo-balanced, i.e. every edge of Tor$[c]$
belongs to some directed cycle. 
\nl
2) Assume that $\bar t \in \text{ pr}({\frak C})$ is balanced,
\ub{then} maj$(\bar t)$ is pseudo-balanced, i.e.
if $t_{x,y} > \frac 12$ then we can find $k \ge 3$ and
$x_0,\dotsc,x_{k-1} \in X$ with no repetitions such
that $(x_0,x_1) = (x,y)$ and $j=i+1$ mod $k \Rightarrow t_{x_i,x_j} >
\frac 12$.
\bigskip

\demo{Proof}  1) As $c$ is weight balanced, it is maj$(\bar t)$ for
some balanced $\bar t \in \text{ pr}({\frak C})$, i.e. $\bar t \in
\text{ pr}^{\text{bl}}({\frak C})$, so now the second conclusion
``every edge of Tor$(c)$ belongs to a directed cycle", follows from
part (2) by the definition of maj$(c)$.   The first conclusion
``partition-balanced" follows from the first; assume $\emptyset
\subsetneqq Y \subsetneqq X,x \in X \backslash Y,y_1 \in Y$ and
$(x,y)$ is an edge of Tor$(c)$ then $t_{u,v} > \frac 12$ hence there
is $k \ge 3$ and $x_0,\dotsc,x_{k-1} \in X$ as in part (2).

Let $i_*$ be the maximal $i \in\{1,\dotsc,k-1\}$ such that $x_i \in
Y$, (well defined as $i=1$ is O.K.) and $j_* = i_* + 1$ mod $k$, so
$t_{x_{i_*},x_{j_*}} > \frac 12,x_{i_*} \in Y,x_{j_*} \in X \backslash
Y$ as required in Definition \scite{b.1}(3A).
\nl
2) Assume $\bar t \in \text{ pr}({\frak C})$ is balanced and $(x,y)$
is an edge of $\bold G_{\text{maj}(\bar t)}$, so 
$t_{x,y} > \frac 12$ and there are no $\langle
x_0,\dotsc,x_{k-1}\rangle$ as promised.  Let $Y_1$ be the set of $z \in X$ such
that there are $k$ and $z_0,\dotsc,z_k \in X \backslash \{y\}$ such
that $z_0 = z,z_k = x,t_{z_i,z_{i+1}} > \frac 12$ for $i<k$.

Let $Y_2 = X \backslash Y_1$ so
\mr
\item "{$\odot_1$}"  $\{x\} \subseteq Y_1 \subseteq X \backslash
\{y\}$ so $(Y_1,Y_2)$ is a partition of $X$ to non-empty sets.
\nl
Now
\sn
\item "{$\odot_2$}"  if $u \in Y_1,v \in Y_2$ then $t_{u,v} \ge \frac
12$.
\nl
[Why?  Toward contradiction assume $t_{u,v} < \frac 12$ so $t_{v,u} >
\frac 12$.  Now by the definition of $u \in Y_1$ there are $z_0 =
u,z_1,\dotsc,z_k = x$ as there; so $\{z_0,\dotsc,z_k\} \subseteq Y_1$
by the definition of $Y_1$ and \wilog \, $\langle
z_0,\dotsc,z_k\rangle$ is without repetitions.  If $v=y$ then 
$y,z_0,z_1,\dotsc,z_k$ is a cycle as required.  If $v \ne y$ then the 
sequence $v,z_0,\dotsc,z_k$ shows that $v \in Y_1$
contradicting the assumption $v \in Y_2$ of $\odot_2$.]
\sn
\item "{$\odot_3$}"  for some $u \in Y_1,v \in Y_2$ we have $t_{u,v} >
\frac 12$.
\nl
[Why?  Choose $u = x,v=y$.]
\ermn
By $\odot_2 + \odot_3$ we get a contradiction to \scite{b.2}(2).
\hfill$\square_{\scite{b.2A}}$
\enddemo
\bigskip

\proclaim{\stag{b.3} Claim}  Assume $|X| \ge 3,{\frak D} \subseteq
{\frak C}^{\text{full}}$ is symmetric, non-empty and balanced.  
\ub{Then}, for any
distinct $x,y,z \in Z$ we have $\bar t^{<x,y,z>} \in 
{\text{\rm pr-cl\/}}({\frak D})$. 
\endproclaim
\bigskip

\demo{Proof}  Let $d \in {\frak D}$ be non-trivial, now Tor$(d) =:
(X,\{(u,v):d\{u,v\} = v\})$ is a directed graph with 
equal out-valance and in-valance
for every node, it has a directed cycle.  As ${\frak D} \subseteq
{\frak C}^{\text{full}}_X$, it follows that this graph has 
a triangle, i.e., $x,y,z \in X$ distinct such that
\mr
\item "{$(*)_1$}"  $d\{x,y\} = y,d\{y,z\} = z,d\{z,x\} = x$.
\ermn
Let $\Pi_{x,y,z} = \{\pi \in \text{ Per}(X):\pi \restriction
\{x,y,z\}$ is the identity$\}$.  Let $\bar t = \Sigma\{\bar t[d^\pi]:
\pi \in \Pi_{x,y,z}\}/|\Pi_{x,y,z}|$. \nl
Clearly $d^\pi \in {\frak D}$ for $\pi \in \Pi_{x,y,z}$ hence $\bar t
\in \text{ pr-cl}({\frak D})$.  Also by $(*)_1$ and the definition
of $\Pi_{x,y,z}$
\mr
\item "{$(*)_2$}"  $t_{x,y} = t_{y,z} = t_{z,x} =1$.
\ermn
Also

$$
\align
|\{w:&w \in X \backslash \{x,y,z\} \text{ and } d\{x,w\} = w\}| = \\
  &|\{w:w \in X \backslash \{x\} \text{ and } d\{x,w\} = w\}| - |\{w:w
\in \{y,z\} \text{ and } d\{x,w\} = w\}| = \\
  &(|X|-1)/2 -1 = (|X| -3)/2 
\endalign
$$
\mn
so

$$
\align
|\{w:w \in X \backslash \{x,y,z\} \text{ and } d\{x,w\} &= x\}| =
 (|X|-3)-|\{w:w \in X \backslash \{x,y,z\} \text{ and} \\
  &d\{x,w\} = w\}| = (|X|-3)-(|X|-3)/2 = (|X|-3)/2. 
\endalign
$$
\mn
hence
\mr
\item "{$(*)_3$}"  $t_{x,w} = 1/2 = t_{w,x}$ 
for $w \in X \backslash \{x,y,z\}$.
\ermn
Similarly
\mr
\item "{$(*)_4$}"  $t_{y,w} = 1/2 = t_{w,y}$ for $w \in X \backslash \{x,y,z\}$
\sn
\item "{$(*)_5$}"  $t_{z,w} = 1/2 = t_{w,z}$ for $w \in X \backslash \{x,y,z\}$
\ermn
and even easier (and as in \S2)
\mr
\item "{$(*)_6$}"  $t_{u,v} = 1/2$ if $u \ne v \in X 
\backslash \{x,y,z\}$.
\ermn
So, by the definition of $\bar t^{<x,y,z>}$, we are done.
\hfill$\square_{\scite{b.3}}$ 
\enddemo
\bigskip

\proclaim{\stag{b.4} Claim}  Assume ${\frak D} \subseteq 
{\frak C}^{\text{full}}$ is symmetric non-empty and 
$c \in {\frak C}$ is pseudo-balanced
\ub{then} $c \in {\text{\rm maj-cl\/}}({\frak D})$.
\endproclaim
\bigskip

\demo{Proof}  Without loss of generality ${\frak D}$ is balanced
(otherwise use \scite{3.1A}).  So by \scite{b.3}
\mr
\item "{$\circledast$}"  if $x,y,z \in X$ are distinct then 
$\bar t^{<x,y,z>} \in \text{ pr-cl}({\frak D})$.
\ermn
Let $\langle \bar x^i:i < i(*) \rangle$
list the set cyc$(c)$ of tuples $\bar x = \langle x_\ell:\ell \le k
\rangle$ such that:
\mr
\item "{$\odot$}"  $(a) \quad k \ge 2,x_\ell \in X$
\sn
\item "{${{}}$}"  $(b) \quad \ell_1 < \ell_2 \le k \Rightarrow x_{\ell_1} \ne
x_{\ell_2}$
\sn
\item "{${{}}$}"  $(c) \quad c\{x_\ell,x_{\ell +1}\} = 
x_{\ell +1} \text{ for } \ell < k$
\sn
\item "{${{}}$}"  $(d) \quad c\{x_k,x_0\} = x_0$.
\ermn
For a tuple $\bar x = \langle x_\ell:\ell < m\rangle$ let $\ell g(\bar
x)$ be the length of $\bar x,m$.

Note
\mr
\item "{$\otimes$}"  for every $\bar x \in \text{ cyc}(c)$ for some
$\bar t = \bar t^{\bar x} \in \text{ pr-cl}({\frak D})$ we have
{\roster
\itemitem{ $(a)$ }  $t_{u,v} = \frac 12 + \frac{1}{2(\ell g(\bar x)-2)} 
\,\,\text{ \ub{if}}$ \nl

\hskip40pt $(u,v) \in \{(x_{\ell_1},x_{\ell_2}):\ell_1 < \ell g(\bar x)-1 \and
\ell_2 = \ell_1 +1 \text{ or}$ \nl

\hskip80pt $\ell_1 = \ell g(\bar x)-1 \and \ell_2=0\}$
\sn
\itemitem{ $(b)$ }  $t_{u,v} = \frac 12 
- \frac{1}{2(\ell g(\bar x)-2)} \,\, \text{ \ub{if} }
(v,u) \text{ is as above}$
\sn
\itemitem{ $(c)$ }  $t_{u,v} = \frac 12 \text{ \ub{if} otherwise}$. 
\endroster}
\ermn
[Why?  If $\bar x = \langle x_\ell:\ell \le k \rangle$, let $\bar t$
be the arithmetic average of \nl
$\langle \bar t^{<x_0,x_1,x_2>},
\bar t^{<x_0,x_2,x_3>},\dotsc,\bar t^{<x_0,x_{k-1},x_k>}\rangle$.] 
\sn
Now let

$$
\bar t = \Sigma\{ \frac{1}{i(*)} \bar t^{\bar x^i}:i < i(*)\}.
$$
\mn
(In fact we just need that $c\{y_0,y_1\} = y_1 \Rightarrow (y_0,y_1)$
appears in at least one cycle $\bar x^i,i < i(*)$).  As every edge of
Tor$(c)$ belongs to a directed cycle easily $c =
\text{ maj}(\bar t)$.  
\nl
${{}}$   \hfill$\square_{\scite{b.4}}$ 
\enddemo
\bn
So now we can give a complete answer.
\demo{\stag{b.5} Conclusion}  Assume
\mr
\item "{$(a)$}"  ${\frak D} \subseteq {\frak C}^{\text{full}}$ is
symmetric, non-empty
\sn
\item "{$(b)$}"  $c \in {\frak C}$.
\ermn
\ub{Then} $c \in \text{ maj-cl}({\frak D})$ \ub{iff} ${\frak D}$ has
a non-balanced member \ub{or} $c$ is pseudo-balanced.
\enddemo
\bigskip

\demo{Proof}  If ${\frak D}$ has no non-balanced member and $c$ is
not pseudo-balanced, by \scite{b.2}(3) we know $c \notin$
maj-cl$({\frak D})$.  \nl
For the other direction, if $d^* \in {\frak D}$ is not balanced 
use \scite{3.1A} that is $(a)' \Leftrightarrow (g)$ 
of claim \scite{2.1} for sym-cl$\{d^*\}$.
Otherwise ${\frak D}$ is balanced non-empty, $c$ is pseudo-balanced and 
we use \scite{b.4}.  \hfill$\square_{\scite{b.5}}$ 
\enddemo
\newpage


\nocite{ignore-this-bibtex-warning} 
\newpage
    
REFERENCES.  
\bibliographystyle{lit-plain}
\bibliography{lista,listb,listx,listf,liste}

\enddocument